\documentclass[a4paper,12pt]{article}
\usepackage{amsfonts,amsmath}
\usepackage{amssymb}
\usepackage[mathscr]{eucal}
\usepackage{amsthm}
\usepackage{amscd}
\usepackage[T2A]{fontenc}
\usepackage{graphics}
\usepackage[dvips]{color}
\usepackage[cp1251]{inputenc}

\newcounter{rem}[section]

\def\text{\mbox}

\def\varkappa{\kappa}
\def\suml{\sum\limits}
\def\maxl{\max\limits}
\def\minl{\min\limits}
\def\intl{\int\limits}
\def\supl{\sup\limits}
\def\infl{\inf\limits}
\def\liml{\lim\limits}
\def\prodl{\prod\limits}

\def\dpfrac{\displaystyle\frac}

\newcommand{\df}{\stackrel{\textrm{def}}{=}}
\newtheorem{thm}{Theorem}
\newtheorem{cor}{Corollary}

\newtheorem{rmk}{Remark}
\def\text{\mbox}

\def\suml{\sum\limits}
\def\maxl{\max\limits}
\def\minl{\min\limits}
\def\intl{\int\limits}
\def\supl{\sup\limits}
\def\infl{\inf\limits}

\def\liml{\lim\limits}
\def\prodl{\prod\limits}
\def\dip{\displaystyle}

\def\dpfrac{\dip\frac}
\title{Retrospective change-point detection and estimation in multivariate linear models}
\date{}
\author{\textbf{Boris Brodsky}\\
Central Institute for Mathematical Economics, RAS, Moscow, Russia
\and\\
\textbf{Boris Darkhovsky}\\
Institute for Systems Analysis RAS, Moscow, Russia}
\date{}
\begin{document}
\maketitle

\begin{abstract}
In this paper the problem of retrospective change-point detection
and estimation in multivariate linear models is considered. The
lower bounds for the error of change-point estimation are proved
in different cases (one change-point: deterministic and stochastic
predictors, multiple change-points). A new method for
retrospective change-point detection and estimation is proposed
and its main performance characteristics (type 1 and type 2
errors, the error of estimation) are studied for dependent
observations in situations of deterministic and stochastic
predictors and unknown change-points. We prove that this method is
asymptotically optimal by the order of convergence of change-point
estimates to their true values as the sample size tends to
infinity. Results of a simulation study of the main performance
characteristics of proposed method in comparison with other well
known methods of retrospective change-point detection and
estimation are presented.

\end{abstract}


\section{Introduction}

This paper deals with change-point problems for multivariate
linear models. We begin with a short review of this field.

The change-point problem for regression models was first
considered by Quandt (1958, 1960). Using econometric examples
Quandt proposed a method for estimation of a change-point in a
sequence of independent observations based upon the likelihood
ratio test.

Let us describe the change-point problem for the linear regression
models considered in the literature. Let $y_1,y_2,\dots,y_n$ be
independent random variables (i.r.v.'s). Under the null hypothesis
$\mathbf{H}_0$ the linear model is
$$
y_i=\mathbf{x}_i^*\mathbf{\beta}+\epsilon_i,\quad 1\le i\le n,
$$
where $\mathbf{\beta}=(\beta_1,\beta_2,\dots,\beta_d)^*$ is an
unknown vector of coefficients, ${\mathbf x}_i^{*}=(1,x_{2i},\\
\dots,x_{di})$ are known predictors (here and below $*$ is the
transposition symbol).

The errors $\epsilon_i$ are supposed to be independent identically
distributed random variables (i.i.d.r.v.'s) with
$\mathbf{E}\epsilon_i=0,\quad 0<\sigma^2=var\,\epsilon_i <
\infty.$

Under the alternative hypothesis $\mathbf{H}_1$ a change at the
instant $k^*$ occurs, i.e.
$$
y_i=\left \{
\begin{array}{rcl}
\mathbf {x}_i^*\mathbf{\beta}+\epsilon_i,& \quad 1\le i\le k^* \\
\mathbf {x}_i^*\mathbf{\gamma}+\epsilon_i,& \quad k^*< i\le n,
\end{array}
\right.
$$
where $k^*$ and $\mathbf{\gamma}\in \mathbb{R}^d$ are unknown
parameters, and $\beta\ne\gamma$.

Denote
$$
\begin{array}{ll}
&\bar {y}_k=\dpfrac{1}{k} \suml_{1\le i\le k}\,y_i,
\bar {\mathbf{x}}_k=\dpfrac{1}{k} \suml_{1\le i\le k}\,\mathbf{x}_i, \\[2mm]
&Q_n=\suml_{1\le i\le n}\,(\mathbf {x}_i-\bar {\mathbf
{x}}_n)(\mathbf{x}_i-\bar {\mathbf{x}}_n)^{*}
\end{array}
$$
and $\mathbf{X}_n=(\mathbf {x}_1,\mathbf {x}_2,\dots,\mathbf
{x}_n)^{*}$,  $Y_n=(y_1,y_2,\dots,y_n)^{*}$.

The least square estimate of $\mathbf{\beta}$ is:
$$
\mathbf {\hat{\beta}_n}=(\mathbf{X}_n^
{*}\mathbf{X}_n)^{-1}\,\mathbf{X}_n^{*} Y_n.
$$

Siegmund with co-authours (James, James, Siegmund (1989)) proposed
to reject $\mathbf{H}_0$ for the large values of $ \maxl_{1\le
k\le n}\,|U_n(k)|, $ where
$$
U_n(k)=(\dpfrac {k}{1-k/n})^{1/2}\,\dpfrac {\bar y_k-\bar
y_n-\mathbf{\hat{\beta}}_n(\mathbf{\bar{
x}_k}-\mathbf{\bar{x}}_n)^{*}} {(1-k(\bar x_k-\bar x_n)(\bar
x_k-\bar x_n)^*/ (Q_n(1-k/n)))^{1/2}}.
$$

Earlier, Brown, Durbin, and Evans (1975) used the cumulative sums
of regression residuals
$$
\suml_{1\le i\le k}\,(y_i-\bar{y}_n-\mathbf{\hat{\beta}}_n
(\mathbf{x}_i-\mathbf{\bar{ x}}_n)^{*}), \quad 1\le k\le n.
$$
It is easy to see that
$$
\begin{array}{ll}
& U_n(k)=w_n(k)\,R_n(k)\\[2mm]
& R_n(k)=(\dpfrac {n}{k(n-k)})^{1/2}\,\suml_{1\le i\le
k}\,(y_i-\bar{y}_n-\mathbf{\hat{\beta}}_n
(\mathbf{x}_i-\mathbf{\bar {x}}_n)^{*})\\[2mm]
&w_n(k)=1-k(\mathbf{\bar{ x}}_k-\mathbf{\bar{ x}}_n)(\mathbf{\bar
{x}}_k-\mathbf{\bar{ x}}_n)^{*}/ (Q_n(1-k/n)))^{-1/2}.
\end{array}
$$

The functionals of $U_n(k)$ and  $R_n(k)$ were used as the test
statistics for detection of change-points in regression
relashionships.

Kim and Siegmund (1989) obtained the limit distribution of
 $\maxl_{1\le k<n}|U_n(k)|$.
Alternatively, Maronna and Yohay (1978), and Worsley (1986) used
the maximum likelihood method for testing $\mathbf{H}_0$ against
$\mathbf{H}_1$ for Gaussian errors. Later Gombay and Horvath
(1994) studied the limit distributions of statistics
$Z_n(i,j)=\maxl_{i\le k<j}|U_n(k)|,\;T_n(i,j)=\maxl_{i\le
k<j}|R_n(k)|$ for deterministic and stochastic regression plans.
The monograph by Csorgo and Horvath (1997) puts together various
results in detection of structural changes in regression models.

Besides change-point detection problems, results in change-point
estimation for regressions are of especial practical importance.
This theme is considered in papers by Darkhovsky (1995), Huskova
(1996), Horvath, Huskova, and Serbinovska (1997). In two last
papers the asymptotical characteristics of change-point estimates
based upon the maximum likelihood statistics are studied. For the
case of contiguous alternatives, the limit distribution of the
change-point estimates is obtained and weak and strong consistency
of these estimates is proved. The paper by Darkhovsky (1995)
develops the nonparametric approach to retrospective change-point
estimation. Here the limit characteristics of change-point
estimates in the functional regression model are studied without
the contiguity assumption, and the rate  of convergence of these
estimates to the 'true' change-point parameters is estimated. Some
generalizations of these results can be found in the monograph by
Brodsky and Darkhovsky (2000).

A new wave of research interest to change-point problems in
regressions was formed in 2000s. Different generalizations to
change-point problems for autoregressive time series (Huskova,
Praskova, Steinebach (2007, 2008), Gombay (2008)), for multiple
change-point estimation in non-stationary time series (Davis, Lee,
Rodriguez-Yam (2006)), for testing change-points in covariance
structure of linear processes (Berkes, Gombay, Horvath (2009))
were studied.

However, as a result we see the multitude of methods proposed for
solving different change-point problems in linear relationships
and almost no theoretical approaches to their \emph{comparative
analysis}. We cannot even estimate the asymptotic efficiency of
these methods. All that is empirically observed for 'structural
breaks' tests in statistics and econometrics can be reduced to the
following 'vague' statement: the power of these methods is rather
low. Let us agree that this 'practical conclusion' requires a more
serious verification.

In this paper, we pursue the following main goals:

1) To prove  the prior theoretical lower bounds for the error
probability in change-point estimation in multivariate models.
These bounds provide the theoretical basis for the proofs of the
asymptotic optimality of change-point estimates and for the
comparative analysis of these estimates;

2) To propose a new nonparametric method for the problem of
retrospective change-point detection and estimation in
multivariate linear systems. Then we study the main performance
characteristics of this method: type 1 and type 2 errors, the
error of change-point estimation.

3) For the problem of multiple change-point detection and
estimation, to propose a general statement in which both {\it the
number of change-points and their coordinates in the sample are
unknown}. For this problem statement, to propose a new
asymptotically optimal method which gives consistent estimates of
an unknown number of change-points and their coordinates.

The structure of this paper is as follows. In Section 2 the
general change-point problem for multivariate linear systems is
formulated and general assumptions are given. In Section 3 we
prove the prior informational inequalities for the main
performance characteristic of the retrospective change-point
problem, namely, the error of change-point estimation. The lower
bounds for the error of estimation are found in different
situations of change-point detection (deterministic and stochastic
regression plan, multiple change-points). In Section 4 we propose
a new method for the retrospective change-point detection and
estimation in multivariate linear models and study its main
performance characteristics (type 1 and type 2 errors, the error
of estimation) in different situations of change-point detection
and estimation (dependent observations, deterministic and
stochastic regression plan, multiple change-points). We prove that
this method is asymptotically optimal by the order of convergence
of change-point estimates to  their true values as the sample size
tends to infinity. In Section 5 a variant of the functional limit
theorem in the case of absence of change-points is given. In
Section 6 a simulation study of characteristics of the proposed
method for finite sample sizes is performed. The main goals of
this study are as follows: to compare performance characteristics
of the proposed method with characteristics of other well known
methods of change-point detection in linear regression models, to
consider more general multivariate linear models and performance
characteristics of the proposed method in these multivariate
models. Section 7 contains main conclusions. All proofs are given
in the Appendix.

\section{Problem statement and general assumptions}

\subsection{General model}

The following basic specification of the multivariate system with
structural changes is considered:
$$
\mathbf{Y}(n)=\mathbf{\Pi} \mathbf{X}(n)+\mathbf{\nu}_n, \quad
n=1,\dots,N   \eqno (1)
$$
where $\mathbf{Y}(n)=(y_{1n},\dots,y_{Mn})^{*}$ is the vector of
endogenous variables, $\mathbf{X}(n)=(x_{1n},\dots,x_{Kn})^{*}$ is
the vector of pre-determined variables, $\Pi$ is $M\times K$
matrix, $\mathbf {\nu}_n=(\nu_{1n},\dots,\nu_{Mn})^{*}$ is the
vector of random errors.

The matrix
$\Pi=\Pi(\vartheta,n),\,\vartheta=(\theta_1,\dots,\theta_k)$ can
change abruptly at some unknown change-points $m_i=[\theta_i
N],\,i=1,\dots,k$ (here and below $[a]$ denote the integer part of
number $a$),  i.e.,
$$
\mathbf{\Pi}(\vartheta,n)=\suml_{i=1}^{k+1}\,\mathbf{a}_i\,\mathbb{I}([\theta_{i-1}N]<n\le
[\theta_i N]),
$$
where $\theta_i$ are unknown change-point parameters such that
$0\equiv \theta_0<\theta_1<\dots \theta_k<\theta_{k+1}\equiv 1$,
$\mathbf{a}_i\ne \mathbf{a}_{i+1},\,i=1,\dots,k$ are unknown
matrices (here and below $\mathbb{I}(A)$ is the indicator of the
set $A$).

The problem is to estimate the unknown parameters $\theta_i$ (and
therefore, the change-points $m_i$) by observations
$\mathbf{Y}(i),\mathbf{X}(i),\,i=1,\dots,N$ (the case
$\theta_i\equiv 1, i=1,\dots,k$ corresponds to the model without
change-points).

Therefore, first, we need to test an obtained dataset of
observations for the presence of change-points. Second, in the
case of a rejected stationarity hypothesis, we wish to estimate
all detected change-points.

 Model (1) generalizes many widely used regression models, namely:

a)\emph{autoregression model (AR)}
$$
y_n=c_0+c_1y_{n-1}+\dots+c_m y_{n-m}+\nu_n,
$$
Here
$\mathbf{X}(n)=(1,y_{n-1},\dots,y_{n-m})^*,\,\mathbf{\Pi}=(c_0,c_1,\dots,c_m)$.

b)\emph{autorgression-moving average (ARMA) model}
$$
y_n=c_1y_{n-1}+\dots+c_k y_{n-k}+d_1 u_{n-\Delta}+\dots+d_m
u_{n-\Delta-m}+\nu_n,
$$
where $u_n$ is the input variable, $y_n$ is the output variable at
the instant $n$, $\Delta$ is the delay time. Here
$\mathbf{X}(n)=(y_{n-1},\dots,y_{n-m},u_{n-\Delta},\dots,u_{n-\Delta-m})^*,\,
\mathbf{\Pi}=(c_1,\dots,c_k, d_1,\dots,d_m)$.

c)\emph{multi-factor regression model}
$$
y_n=c_1y_{n-1}+\dots+c_k
y_{n-m}+\suml_{i=1}^r\,\suml_{j=1}^{l_i}\,d_{ij}x_i(n-j)+\nu_n,
$$
where $r,m,l_i \ge 1$. Here
$\mathbf{X}(n)=(y_{n-1},\dots,y_{n-m},x_1(n-1),\dots,x_1(n-l_1),x_2(n-1),
\dots,x_2(n-l_2),\dots,
x_r(n-1),\dots,x_r(n-l_r))^*,\,\mathbf{\Pi}=(c_1,\dots,c_k,d_{11},\\
\dots,d_{rl_r})$.

d)\emph{simultaneous equation systems (SES)}
$$
B\mathbf{Y}(n)+\Gamma \mathbf{X}(n)=\mathbf{\epsilon}_n,
$$
where $\mathbf{Y}(n)=(y_{1n},y_{2n},\dots,y_{Mn})^{*}$ is the
vector of endogenous variables,
$\mathbf{X}(n)=(x_{1n},x_{2n},\dots,x_{Kn})^{*}$ is the vector of
pre-determined variables (all exogenous variables plus lagged
endogenous variables),
$\epsilon_n=(\epsilon_{1n},\epsilon_{2n},\dots,\epsilon_{Mn})^{*}$
is the vector of random errors, $B$ is a $M\times M$
non-degenerate matrix ($\det B\ne 0$), $\Gamma$ is a $M\times K$
matrix.

This general structural form of the SES can be written in the
following {\it reduced form}:
$$
\mathbf{Y}(n)=-B^{-1}\,\Gamma
\mathbf{X}(n)+B^{-1}\epsilon_n=\mathbf{\Pi}\mathbf{X}(n)+\nu_n
$$

This system is usually used for the analysis of change-points
(structural changes) in multivariate linear models (see, e.g.,
Bai, Lumsdaine, Stock (1998)).

\subsection{General assumptions}

In this subsection we formulate general assumptions which will be
used in our main theorems 3-5. Some specific assumptions will be
formulated together with the corresponding theorems.

Let us start from the following definitions. Consider the
probability space ($\Omega,\mathfrak{F},\mathbf{P}$). Let
$\mathcal{H}_{1}$ and $\mathcal{H}_{2}$ be two $\sigma $-algebras
from $\mathfrak{F}$. Consider the following measure of dependence
between $\mathcal{H}_{1}$ and $\mathcal{H}_{2}$:
$$
 \psi (\mathcal{H}_{1}, \mathcal{H}_{2})  = \sup_{A \in
\mathcal{H}_{1},B \in \mathcal{H}_{2}, \mathbf{P}(A)
\mathbf{P}(B)\neq 0} \Big\vert \frac{\mathbf{P}(AB)}{\mathbf{P}(A)
\mathbf{P}(B)}-1 \Big\vert
$$

Suppose ($X_{i},i\ge 1$) is a sequence of random vectors defined
on ($\Omega,\mathfrak{F},\mathbf{P}$). Denote by
$\mathfrak{F}^{t}_{s}=\sigma \{X_{i}: s\le i\le t\}, 1\le s\le t<
\infty$ the minimal $\sigma $-algebra generated by random vectors
$X_{i}, s \le i \le t$. Define
$$
 \psi (n)  = \sup_{t\ge 1} \psi (\mathfrak{F}^{t}_{1},
\mathfrak{F}^{\infty }_{t+n})
$$

A) \emph{Mixing condition}

We say that scalar random sequence $\{x_n\}$ satisfies the $\psi
$-\emph{mixing condition} if the function $\psi(n)$ (which is also
called the
 \emph{$\psi $-mixing coefficient}) tends to zero as $n$ goes to infinity.

We say that vector random sequence
$\{X(n)\},\,X(n)=\left(x_1(n),\dots,x_k(n)\right)^*$ satisfies the
\emph{uniform $\psi$-mixing condition} if
$\maxl_{i,j}\psi_{ij}(n)$ tends to zero as $n$ goes to infinity,
where $\psi_{ij}(n)$ is the $\psi$-mixing coefficient for the
sequence $\{x_i(n)x_j(n)\}$.

The $\psi$-mixing condition is satisfied in most practical
situations of change-point detection. In particular, for a Markov
chain (not necessarily stationary), if $\psi(n)<1$ for a certain
$n$, then $\psi(k)$ goes to zero at least exponentially as
$k\to\infty$ (see Bradley, 2005, theorem 3.3).

B) \emph{Cramer condition}

Let $\{\zeta(n)\},\,\,\zeta(n)=
\left(\zeta_{1}(n),\dots,\zeta_{k}(n)\right)^*$ be a vector random
sequence. We say that the \emph{uniform Cramer condition} is
satisfied if there exists a constant $L>0$ such that
$$
\supl_n\,\mathbf{E}\exp \left(t \zeta_i(n)\zeta_j(n)\right)<
\infty
$$
for every $i,j=1,\dots,k$ and $|t|< L$.

For a centered random sequence $\xi_n$ this condition is
equivalent to the following: there exist constants $g>0,\,T>0$
such that for each $|t|<T$:
$$
\supl_n\mathbf{E}e^{t\xi_n}\le \exp\left(\frac{t^2g}{2}\right).
$$

\section{Preliminary results: prior inequalities}

\subsection{Unique change-point}

On a probability space $(\Omega,{\cal F},\mathbf{P}_{\theta})$
consider a sequence of i.r.v.'s $x_1,\dots,x_{\scriptscriptstyle
N}$ with the following density function (w.r.t. some
$\sigma$-finite measure $\mu$)
$$
f(x_n)=\left \{
\begin{array}{ll}
& f_0(x_n,n/N), \qquad 1\le n\le [\theta N], \\
& f_1(x_n,n/N), \qquad [\theta N] < n\le N.
\end{array}
\right.  \eqno (2)
$$
Here $0< \theta <1$ is an \emph{unknown change-point parameter}.

Define the following objects:
$$
T_{\scriptscriptstyle N}(\Delta):\mathbb{R}^{\scriptscriptstyle
N}\longrightarrow \Delta\subset\mathbb{R}^1 \eqno(3)
$$
is the Borel function on $\mathbb{R}^{\scriptscriptstyle N}$ with
the values in the set $\Delta$;
$$
\mathcal{M}_{\scriptscriptstyle N}(\Delta)=\{T_{\scriptscriptstyle
N}(\Delta)\} \eqno(4)
$$
is the collection of all Borel functions $T_{\scriptscriptstyle
N}$.

\begin{thm}
 Suppose the following assumption is satisfied:

the functions
$J_0(t)\df\mathbf{E}_0\ln\dpfrac{f_0(x,t)}{f_1(x,t)}$ and
$J_1(t)\df\mathbf{E}_1\ln\dpfrac{f_1(x,t)}{f_0(x,t)}$ are
continuous at $[0,1]$ and such that
$$
J_0(t)\ge\delta>0,\,\,J_1(t)\ge\delta>0.
$$

Then for any fixed $0<\theta <1,\,0< \epsilon < \theta \wedge
(1-\theta)$ the following inequality holds:
$$
\liminf_{N\to\infty}\,N^{-1}\ln\infl_{\hat\theta_N\in
{\mathcal{M}}_N((0,1))}\,
\mathbf{P}_{\theta}\{|\hat\theta_N-\theta|> \epsilon \} \ge
-\min\left
(\intl_{\theta}^{\theta+\epsilon}\,J_0(t)dt,\,\intl_{\theta-
\epsilon}^{\theta}\,J_1(t)dt\right).
$$
\end{thm}

The proof of this theorem is given in the Appendix A.

\begin{rmk}
 The lower bound in Theorem 1 can not be improved essentially. It
follows from the results of Korostelev (1997). In this work the
exact lower bound for the change-point estimate in continuous time
model for the Wiener process was given. The exact lower bound in
Korostelev (1997) differs from our bound only by a constant
factor.
\end{rmk}

\bigskip
Consider the following particular cases of model (2).

\medskip
\emph{1. A break in the trend function $\phi(t)$ of the
mathematical expectation  of Gaussian observations}

Let
$$\begin{array}{rcl}
f_0(x,t)&=&h(x)\exp \left(\phi_0(t) x- \phi_0^2(t)/2\right),\quad
t\le \theta
\\[2mm]
f_1(x,t)&=&h(x)\exp \left(\phi_1(t) x-\phi_1^2(t)/2\right),\quad
t> \theta,
\end{array}
$$
where $h(x)=\dpfrac {1}{\sqrt {2\pi}}\exp (-x^2
/2),\,\phi_0(\cdot)\not=\phi_1(\cdot)$.

In this case from Theorem 1 we obtain the following lower bound
for the error probability:
$$\begin{array}{ll}
&\mathbf{P}_{\theta}\{|\hat\theta_N-\theta|  >  \epsilon \} \ge
(1-o(1))\cdot\\
&\cdot\exp \left(-\dpfrac N2 \,\min
\big(\intl_{\theta}^{\theta+\epsilon}\,(\phi_0(t)-\phi_1(t))^2 dt,
\intl_{\theta-\epsilon}^{\theta}\,(\phi_0(t)-\phi_1(t))^2\,dt\big)\right).
\end{array}
$$

\bigskip
{\it 2. Linear regression with deterministic predictors and
Gaussian errors}

Let
$$
y_n=c_1(n)x_{1n}+\dots+c_k(n)x_{kn}+\xi_n, \quad n=1,\dots,N,
\eqno(5)
$$
where $\{\xi_{n}\}$ is a sequence of independent Gaussian r.v.'s
with zero mean,
 $\xi_n\sim {\mathcal{N}}(0,\sigma^2),\,\, \\ \mathbf
{c}(n)\df(c_1(n),\dots,c_k(n))^{*}=\mathbf {a}\mathbb{I}(n\le [\theta N])+\mathbf {b}\mathbb{I}(n>[\theta N]),
\,\,\mathbf{a}=(a_1,\dots,a_k)^*\not=\linebreak\mathbf{b}=(b_1,\dots,b_k)^*,\,\, x_{in}=f_i( n/N),\,n=1,\dots, N$, and $f_i(\cdot)\in
C[0,1],\,i=1,\dots,k$.

In this case from Theorem 1 applied to the sequence of
observations $y_1,\dots,y_N$ we obtain:
$$\begin{array}{ll}
&\mathbf{P}_{\theta}\{|\hat\theta_N - \theta|> \epsilon \}\ge (1-o(1))\cdot \\
& \cdot\exp\left(-\dpfrac N{2\sigma^2}\, \min
\big(\intl_{\theta}^{\theta+\epsilon}\,(\suml_{i=1}^k\,f_i(t)(a_i-b_i))^2
dt,\intl_{\theta-\epsilon}^{\theta}\,(\suml_{i=1}^k\,f_i(t)(a_i-b_i))^2
dt\big) \right).
\end{array}
$$

\bigskip
{\it 3. Linear stochastic regression model with Gaussian
predictors}

Consider model (5) with $\xi_{n}\equiv 0$. Suppose that there exist continuous functions $f_i(\cdot), \sigma_i(\cdot), \,i=1,\dots, k$ such that
$x_{in}$ are Gaussian i.r.v.'s, $x_{in}\sim {\mathcal {N}}\left(f_i(n/N),\sigma_i^2(n/N)\right),\linebreak n=1,\dots,N$. Suppose also that
$x_{in}$ and $x_{jn}$ are independent for $i\not=j$ and $\mathbf {c}(n)$ is the same as in model (5).

Then from Theorem 1 we obtain:
$$
\mathbf{P}_{\theta}\{|\hat\theta_N-\theta|> \epsilon\}\ge
(1-o(1))\exp \left(-\dpfrac N2 \min
\big(\intl_{\theta}^{\theta+\epsilon}\,J_0(t)dt,
\intl_{\theta-\epsilon}^{\theta}\,J_1(t)dt\big)\right),
$$
where
$$
\begin{array}{ll}
&J_0(t)=\left(\dpfrac {\phi_0(t)}{\Delta_0(t)}-\dpfrac
{\phi_1(t)}{\Delta_1(t)}\right)^2+ 2\dpfrac
{\phi_0(t)}{\Delta_0(t)}\dpfrac {\phi_1(t)}{\Delta_1(t)}\left(1-
\dpfrac {\Delta_0(t)}{\Delta_1(t)}\right)+\\[3mm]
&2\ln\dpfrac {\Delta_1(t)}{\Delta_0(t)}+ \left(1+\dpfrac
{\phi_0^2(t)}{\Delta_0^2(t)}\right)\left(\dpfrac
{\Delta_0(t)}{\Delta_1(t)}-1\right),
\end{array}
$$
and
$$\begin{array}{ll}
&\phi_0(t)=a_1f_1(t)+\dots+a_kf_k(t),\,\Delta_0^2(t)=a_1^2\sigma_1^2(t)+
\dots+a_k^2\sigma_k^2(t),\\
&\phi_1(t)=b_1f_1(t)+\dots+b_kf_k(t),\,\Delta_1^2(t)=b_1^2\sigma_1^2(t)+
\dots+b_k^2\sigma_k^2(t).
\end{array}
$$

\subsection{Multiple change-points}

Theorem 1 can be generalized to the case of several change-points
in the sequence of independent r.v.'s with the following density
function:
$$
f(x_n)=f_i(x_n,n/N)\,\mathbb{I}([\theta_{i-1}N]<n\le [\theta_i
N]),\quad n=1,\dots,N,
$$
where $i=1,\dots,k+1$ and $0\equiv \theta_0< \theta_1< \dots
<\theta_k <\theta_{k+1}\equiv 1$.

Suppose the following assumptions are satisfied:

i) change-points $\theta_i$ are such that $\minl_{1\le i\le
k+1}(\theta_i-\theta_{i-1})\ge \delta >0$.

ii) the functions $J_i(t)=\mathbf{E}_i\,\ln\dpfrac
{f_i(x,t)}{f_{i-1}(x,t)}$ and
$J^{i-1}(t)=\mathbf{E}_{i-1}\,\ln\dpfrac
{f_{i-1}(x,t)}{f_i(x,t)}$, $i=1,\dots,k$ are continuous at $[0,1]$
and such that
$$
J_i(t)\ge\Delta>0,\,i=1,\dots,k
$$

For the multiple change-point problem we estimate both the number
$k$ and the vector $\vartheta\df (\theta_1,\dots,\theta_k)$ of
change-points' coordinates. Let $s^*\df[1/\delta]$ and denote
$Q=\{1,2,\dots,s^*\}$.

For any $s\in Q$ define
$$
\begin{array}{ll}
&\mathcal{D}_{s} = \{x\in \mathbb{R}^{s}: \delta \le x_{i}\le
1-\delta,\,
 x_{i+1} - x_{i} \ge \delta , x_{0}\equiv 0, x_{s+1} \equiv  1\}\\
&\mathcal{D}^{\star} =\bigcup^{s^{\star}}_{i=1}\mathcal{D}_{i},
\mathcal{D}^{\star} \subset\mathbb{R}^{s^{\star}}\equiv
\mathbb{R}^{\star}
\end{array}
\eqno(6)
$$
By the construction, an unknown vector $\vartheta$ is an arbitrary
point of the set $\mathcal{D}_{k}$ and an unknown number of the
change-points $k$ is an arbitrary point of the set $Q$.

As before, it is reasonable to consider objects (3)-(4). In this
notation $\mathcal{M}_{\scriptscriptstyle N}(\mathcal{D}^*)$ is
the set of all arbitrary estimates of the parameter $\vartheta$
and $\mathcal{M}_{\scriptscriptstyle N}(Q)$ is the set of all
arbitrary estimates of the parameter $k$ on the basis of
observations with the sample size $N$.

Let $\hat{k}\in\mathcal{M}_{\scriptscriptstyle N}(Q)$ is an
estimate of an unknown number of change-points $k$ and
$\hat{\vartheta}\in\mathcal{M}_{\scriptscriptstyle
N}(\mathcal{D}_k)$ is an estimate of unknown change-point
coordinates on condition that the number of the coordinates was
estimated correctly.

\begin{thm}
 Suppose assumptions i) and ii) are satisfied. Then for any fixed
$0<\epsilon<\delta$ the following inequality holds:
$$\begin{array}{ll}
& \liminf_{N\to\infty} N^{-1} \ln
\infl_{\hat{\vartheta}\in\mathcal{M}_{\scriptscriptstyle
N}(\mathcal{D}_k)}\infl_{\hat{k}\in\mathcal{M}_{\scriptscriptstyle
N}(Q)}\supl_{\vartheta\in\mathcal{D}_k}\supl_{k\in Q}
\mathbf{P}_{\theta}\{\{\hat k\ne k\}\cup \{(\hat k=k)\cap \\
& \cap (\maxl_{1\le i\le k}|\hat\theta_i-\theta_i|
> \epsilon )\}  \ge - \minl_{1\le i\le k}\min
(\intl_{\theta_i}^{\theta_i+\epsilon}\,J^{i-1}(\tau)d\tau,\;
\intl_{\theta_i-\epsilon}^{\theta_i}\,J_i(\tau)d\tau).
\end{array}
$$
\end{thm}

The proof of this theorem is given in the Appendix B.

\section{Main results}
Now consider model (1). In this Section we assume that the uniform
mixing condition (A) and the uniform Cramer condition (B) (see
Section 2) are satisfied, and an unknown vector of change-point
parameters $\vartheta=(\theta_1,\dots,\theta_k)$ is such that
$0<\beta\le\theta_1<\theta_2<\dots<\theta_k\le\alpha<1$, where
$\beta,\,\alpha$ are known numbers. Everywhere below the measure
$\mathbf{P}_{\vartheta}$ corresponds to a sample with the
change-point $\vartheta$ ($\mathbf{P}_0$ corresponds to a sample
without change-points).
\subsection{Unique change-point}
In this subsection model (1) with unique change-point
$0<\beta\le\theta\le\alpha<1$ is considered.
\subsubsection{Deterministic predictors}
Let us formulate  assumptions for model (1) in the case of a
unique change-point (remind that in model (1) the vector
$\mathbf{X}(n)$ has the dimension $K$ and the vector
$\mathbf{Y}(n)$ has the dimension $M$):

a) the vector random sequence $\{\nu_n\}$ satisfies conditions (A)
and (B) (see section 2).

b) there exist functions $f_i(\cdot)\in C[0,1],\,i=1,\dots,K$ such
that $x_{in}=f_i(n/N),n=1,\dots,N$.

Denote $F(t)=\left(f_1(t),\dots,f_K(t)\right)^{*},\,t\in[0,1]$.

c) for arbitrary $0\le t_1<t_2\le 1$, the matrix
$$
A(t_1,t_2)\df\displaystyle\int_{t_{1}}^{t_{2}}F(s)F^{*}(s)ds
$$
is positive definite (below we denote $A(t)\df A(0,t),\,\,A(1)\df
I$).

In virtue of our assumptions, the matrix $I$  is symmetric and
positive definite.

Define  $K\times M$ matrix
$$
Z(n_1,n_2)=\suml_{i=n_1}^{n_2}\,F(i/N)\mathbf{Y}^*(i)
$$
and  $K\times K$ matrix
$$
\mathcal{P}_{n_{1}}^{n_{2}}\df
\suml_{k=n_1}^{n_2}\,F(k/N)F^*(k/N),\quad    1\le n_1<n_2\le N.
$$

The following matrix statistic is used for estimation of an
unknown change-point:
$$
\mathcal{Z}_N(n)=N^{-1}\left( Z(1,n)-\mathcal{P}_{1}^{n}
(\mathcal{P}_1^{N})^{-1}\,Z(1,N)\right).  \eqno (7)
$$

An arbitrary point $\hat n$ of the set $arg\maxl_{[\beta N]\le
n\le [\alpha N]}\,\|\mathcal{Z}_{\scriptscriptstyle N}(n)\|^2$ is
assumed to be the estimate of an unknown change-point (here and
below $\| C\|$ denotes the Gilbert norm of a quadratic matrix $C$,
namely $\|C\|=\sqrt{tr(CC^*)}$).

We define also the value $\hat\theta_N=\hat {n}/N$ - the estimate
of the change-point parameter $\theta$.

Denote $B\df
B(\theta)=\left(E-I^{-1}A(\theta)\right)(\mathbf{a}-\mathbf{b})^*$.

\begin{thm}
Suppose assumptions a)--c) are satisfied and $\textrm{rank}(B)=M$
if $\theta\in [\beta,\alpha]$.

Then the estimate $\hat\theta_N$ converges to the change-point
parameter $\theta \,\,\mathbf{P}_{\theta}$-almost surely as
$N\to\infty$.

Besides, for any fixed $(\alpha-\beta)>\epsilon>0$ the following
inequality is satisfied for $N>N_0(F)$:
$$
\supl_{\beta\le\theta\le\alpha}\mathbf{P}_{\theta}\{|\hat\theta_{\scriptscriptstyle
N}-\theta|> \epsilon \} \le
m_0\left(C(\epsilon,N)/\mathcal{R}\right)\left \{
\begin{array}{ll}
& \exp \left(-\dpfrac {N\beta
\Big(C(\epsilon,N)/\mathcal{R}\Big)^2}
{4gm_0\left(C(\epsilon,N)/\mathcal{R}\right)}\right),\\
& \textrm{  if  } C(\epsilon,N) \le \mathcal{R}gT \\
& \exp\left(-\dpfrac
{TN\beta\Big(C(\epsilon,N)/\mathcal{R}\Big)}{4m_0\left(C(\epsilon,N)/\mathcal{R}\right)}\right),
\\
&\textrm{ if } C(\epsilon,N)>\mathcal{R}gT.
\end{array}
\right. \eqno(8)
$$
where the constants $g,\,T,\,m_0(\cdot)\ge 1$ are taken from the
uniform Cramer's and $\psi$-mixing conditions, respectively,
$C(\epsilon,N)=[\dpfrac{\epsilon\lambda_{\scriptscriptstyle
F}}{4\mathcal{M}} \|\mathbf
{a}-\mathbf{b}\|^2-L_{\scriptscriptstyle_{F}}/N\Big],\,\,\,
N_0(F),\,\lambda_{\scriptscriptstyle F},\,L_{\scriptscriptstyle
F},\,\mathcal{R}$ are constants which can be exactly calculated
for any given family of functions $F(t)$, and the constant
$\mathcal{M}$ is given in the proof.
\end{thm}
\begin{rmk}
The assumption $\textrm{rank}B=M$ yields $K\ge M$, i.e., the
number $M$ of endogenous variables in (1) cannot exceed the number
$K$ of pre-determined variables. Note that for one regression
equation this assumption is always satisfied.
\end{rmk}
\begin{rmk}
For independent random errors $m_0(\epsilon)=1$.
\end{rmk}
\begin{rmk}
Comparing theorems 1 and 3, we conclude that the order of
convergence of the proposed estimate of the change-point parameter
to its true value is asymptotically optimal as $N\to\infty$.
\end{rmk}
\begin{rmk}
For any given family of functions $F(t)$ one can calculate the
function
$f(t)=\|m(t)\|^2,\,m(t)=\liml_{N\to\infty}\mathbf{E}_\theta\mathcal{Z}_{\scriptscriptstyle
N}([Nt])$ (see the proof) and investigate this function on the
square $(\theta,t)\in[\beta,\alpha]\times [\beta,\alpha]$. Such
investigation gives the opportunity to calculate all constants
from the formulation.
\end{rmk}

The proof of Theorem 3 is given in the Appendix C.

\medskip
From the proof we obtain the following
\begin{cor}
Let $C>0$ be the decision threshold and $\mathbb{C}\df
C-\dpfrac{L_{\scriptscriptstyle F}}{N}$. Then:

- for type 1 error the following inequality is satisfied:
$$
\mathbf{P}_0 \{\maxl_{[\beta N]\le n\le [\alpha N]}
 \|{\cal Z}_{\scriptscriptstyle N}(n)\|^2>C \}\le m_0\left(\mathbb{C}/\mathcal{R}\right)
\left \{
\begin{array}{ll}
& \exp \left(-\dpfrac
{TN\mathbb{C}\beta}{4\mathcal{R}m_0\left(\mathbb{C}/\mathcal{R}\right)}\right), \\
& \textrm{ if } \mathbb{C} >\mathcal{R}gT \\
& \exp \left(-\dpfrac {N\beta \mathbb{C}^2}{4\mathcal{R}^2g
m_0\left(\mathbb{C}/\mathcal{R}\right))}\right), \\
& \textrm{ if } \mathbb{C} \le \mathcal{R}gT,
\end{array}
\right. \eqno(9)
$$

- for type 2 error  the following inequality is satisfied:
$$
\mathbf{P}_{\theta}\{\maxl_{[\beta  N]\le n\le [\alpha N]} \|{\cal
Z}_{\scriptscriptstyle N}(n)\|^2\le C \} \le m_0(d) \left \{
\begin{array}{ll}
& \exp \left(-\dpfrac {TN\beta d}{4m_0(d)}\right), \ \ d >gT \\
& \exp\left(-\dpfrac {N\beta d^2}{4g m_0(d)}\right), \ d \le gT,
\end{array}
\right.
$$
where
$d=\mathcal{R}^{-1}\left(\|m(\theta)\|-C-\dpfrac{L_{\scriptscriptstyle
F}}{N}\right)>
0,\,\,\|m(\theta)\|^2=\textrm{tr}(B^{*}A^2(\theta)B)$.
\end{cor}
\subsubsection{Stochastic predictors}

In this subsection we suppose that predictors $x_{ji}$ in (1) are
random. On the probability space $(\Omega,\mathcal
{F},\mathbf{P}_{\theta})$ consider filtration $\{{\mathcal
{F}}_n\},\,n=1,\dots,n$, where $\{{\mathcal {F}}_n\}\in\mathcal
{F},\,\, {\mathcal {F}}_n$ can be interpreted as all available
information up to the instant $n$.

Put $\mathbf{X}(n)\df \left(x_{1n},\dots,x_{Kn}\right)^*$.

Suppose that the following conditions are satisfied:

a) there exists a continuous symmetric matrix function $V(t), t\in
[0,1]$ such that the matrix $\intl_{t_1}^{t_2} V(s)ds$ is positive
definite for any $0\le t_1<t_2\le 1$, and
$\mathbf{E}_{\theta}\mathbf{X}(n)\mathbf{X}^*(n)=V(n/N)$;

b) the sequence of random vectors
$\{\left(\mathbf{X}(n),\nu_{n}\right)\}$ satisfies the uniform
Cramer's and $\psi$-mixing conditions;

c) the random sequence $\{\nu_{n}\}$ is a martingale-difference
sequence w.r.t. the filtration  $\{\mathcal{F}_{n}\}$;

d) the vector of predictors $\mathbf{X}(n)\df
\left(x_{1n},\dots,x_{Kn}\right)^*$ is
$\mathcal{F}_{n-1}$-measurable.

On the segment $[0,1]$ define the $K\times M$ matrix process
$$
u_{\scriptscriptstyle
N}(t)\df\suml_{i=1}^{[Nt]}\,\mathbf{X}(i)\mathbf{Y}^*(i),
$$
and the $K\times K$ matrix process
$$
\mathcal{T}_{\scriptscriptstyle N}(t)\df
\suml_{k=1}^{[Nt]}\,\mathbf{X}(k)\mathbf{X}^*(k).
$$

In virtue of conditions a), b), c), the matrix process
$N^{-1}\mathcal{T}_{\scriptscriptstyle N}(t)$ weakly converges (in
the Skorokhod space) to a positive definite symmetric matrix
function $\mathbb{R}(t)\df\displaystyle\int_{0}^{t}V(s)ds$, and
the rate of convergence is exponential. Below we denote
$\mathbb{R}(1)\df\mathbb{R}$.

Analogously, due to conditions a)-d), the matrix process
$N^{-1}\suml_{k=1}^{[Nt]}\mathbf{X}(k)\nu^*(k)$ weakly converges
to zero with the exponential rate. Both conclusions follow from
the fact that the random processes
$$\begin{array}{rcl}
N^{-1}\suml_{n=1}^{[Nt]}
\left(x_{in}x_{jn}-\mathbf{E}_{\theta}x_{in}x_{jn}\right), \\
N^{-1}\suml_{n=1}^{[Nt]}
\left(x_{in}\nu_{n}\right),\,\,i,j=1,\dots,k
\end{array}
$$
weakly converge to zero (as $N\to\infty$) with the exponential
rate (see Brodsky, Darkhovsky (2000)).

For estimation of an unknown change-point, the following statistic
is used:
$$
\mathbb{Z}_{\scriptscriptstyle N}(n)=N^{-1}\Big(
u_{\scriptscriptstyle N}(n/N)-\mathcal{T}_{\scriptscriptstyle
N}(n/N) (\mathcal{T}_{\scriptscriptstyle
N}(1))^{-1}\,u_{\scriptscriptstyle N}(1)\Big),\,n=1,2,\dots,N.
\eqno(10)
$$

An arbitrary point  $\hat n$ of the set $\mathrm{Arg}\maxl_{[\beta
N]\le n\le [\alpha N]}\,\|\mathbb{ Z}_{\scriptscriptstyle
N}(n)\|^2$ is assumed to be the estimate of an unknown
change-point. Again we define $\hat\theta_{\scriptscriptstyle
N}=\hat {n}/N$ as the estimate of the change-point parameter
$\theta$.

Statistic (10) generalizes statistic (7) to the situation of
stochastic predictors. Assumptions a)-d) guarantee the analogous
properties of this statistic. In particular, the limit value (as
$N\to\infty$) of the mathematical expectation of the statistic
$\mathbb{Z}_{\scriptscriptstyle N}([Nt])$ attains its unique
global maximum on the segment $[0,1]$ at the point $t^*=\theta$.

Assumptions a)-d) guarantee convergence in probability of an
arbitrary point of $\mathrm{Arg}\maxl_{[\beta N]\le n\le [\alpha
N]}\|\mathbb{Z}_{\scriptscriptstyle N}(n)\|^2$ to the point
$\theta$ with the exponential rate. Hence  the
$\mathbf{P}_{\theta}$-a.s. convergence of the proposed estimate to
$\theta$ follows.

\begin{thm}
Suppose that the conditions a)-d) are satisfied and
$\mathrm{rank}(\mathbb{B})=M$ if $\theta\in [\beta,\alpha]$, where
$\mathbb{B}\df\mathbb{B}(\theta)=\Big(E-\mathbb{R}^{-1}\mathbb{R}(\theta)\Big)(\mathbf{a}-\mathbf{b})^*$.

Then the estimate $\hat\theta_{\scriptscriptstyle N}$ of the
change-point parameter $\theta$ converges to
$\theta\,\,\mathbf{P}_{\theta}$-a.s. as $N\to\infty$.

Besides, there exists the number $N_1=N_1(\{\mathbf{X}(n)\})$ such
that for $N>N_1$ and any fixed
$\epsilon,\,\,(\min\left((\alpha-\beta),\|\mathbb{R}\|/2\right)>\epsilon>0)$,
the following inequality holds:
$$
\begin{array}{ll}
&\supl_{\beta\le\theta\le\alpha}\mathbf{P}_{\theta}\{|\hat\theta_{\scriptscriptstyle
N}-\theta|> \epsilon
\} \le \delta_{\scriptscriptstyle N}(\epsilon)+\\[3mm]
&m_0\left(\mathbb{C}(\epsilon,N)/\mathbf{R}\right)\left \{
\begin{array}{ll}
& \exp \left(-\dpfrac {N\beta
\Big(\mathbb{C}(\epsilon,N)/\mathbf{R}\Big)^2}
{4gm_0\left(\mathbb{C}(\epsilon,N)/\mathbf{R}\right)}\right), \textrm{  if  } \mathbb{C}(\epsilon,N) \le \mathbf{R}gT \\
& \exp\left(-\dpfrac
{TN\beta\Big(\mathbb{C}(\epsilon,N)/\mathbf{R}\Big)}{4m_0\left(\mathbb{C}(\epsilon,N)/\mathbf{R}\right)}\right),
\textrm{ if } \mathbb{C}(\epsilon,N)>\mathbf{R}gT,
\end{array}
\right.
\end{array}
$$
where
$\mathbb{C}(\epsilon,N)=\Big[\dpfrac{\epsilon\lambda_{\scriptscriptstyle
V}}{4\mathbb{M}}
\|\mathbf{a}-\mathbf{b}\|^2-\dpfrac{L_{\scriptscriptstyle
V}}{N}\Big],\,\mathbb{M}=\maxl_{\beta\le t\le\alpha}\|M(t)\|$, the
constants $g,T,m_0(\cdot)$ are taken from the uniform Cramer's and
$\psi$-mixing conditions, and $M(t),\,\lambda_{\scriptscriptstyle
V},\,L_{\scriptscriptstyle V},\, \delta_{\scriptscriptstyle
N},\,\mathbf{R}$ are described in the proof.

In particular, for independent observations $m_0(\cdot)=1$.
\end{thm}

Comparing Theorems 1 and 3, we conclude that the order of
convergence of the proposed estimate of the change-point parameter
to its true value is asymptotically optimal as $N\to\infty$.

The proof of Theorem 4 is given in the Appendix D.

\medskip
From the proof we obtain the following
\begin{cor}
Let $S>0$ be the decision threshold and $\mathbb{S}\df
S-\dpfrac{L_{\scriptscriptstyle V}}{N}$. Then:

- for type 1 error the following inequality is satisfied:
$$
\mathbf{P}_0 \{\maxl_{[\beta N]\le n\le [\alpha N]}
\|\mathbb{Z}_{\scriptscriptstyle N}(n)\|^2>S \}\le
\delta_{\scriptscriptstyle
N}(\mathbb{S})+m_0\left(\mathbb{S}/\mathbf{R}\right) \left \{
\begin{array}{ll}
& \exp \left(-\dpfrac {TN\mathbb{S}\beta}{4\mathbf{R}m_0\left(\mathbb{S}/\mathbf{R}\right)}\right), \\
&\mathbb{S} >\mathbf{R}gT \\
& \exp \left(-\dpfrac {N\beta \mathbb{S}^2}{4\mathbf{R}^2g m_0\left(\mathbb{S}/\mathbf{R}\right)}\right), \\
&\mathbb{S} \le \mathbf{R}gT.
\end{array}
\right.
$$

- for type 2 error the following inequality holds:
$$
\mathbf{P}_{\theta}\{\maxl_{[\beta N]\le n\le [\alpha N]}
\|\mathbb{Z}_{\scriptscriptstyle N}(n)\|^2\le S \}
\le\delta_{\scriptscriptstyle N}(\mathbb{S})+ m_0(r) \left \{
\begin{array}{ll}
& \exp \left(-\dpfrac {TN\beta r}{4Rm_0(r)}\right), \\
&r >RgT \\
& \exp\left(-\dpfrac {N\beta r^2}{4R^2g m_0(d)}\right), \\
&r \le RgT,
\end{array}
\right.
$$
where
$r=\mathbf{R}^{-1}\left(\|M(\theta)\|-S-L_{\scriptscriptstyle
V}\right)>0;\;\|M(\theta)\|^2=\mathrm{tr}(\mathbb{B}^*\mathbb{R}^2(\theta)\mathbb{B})$.
\end{cor}

\subsection{Multiple change-points}

The proposed method can be generalized to problems of detection
and estimation of multiple change-points in regression models. A
widespread approach to solving these problems (see, e.g.,  Bai,
Lumsdaine, Stock (1998)) consists in decomposition of the whole
obtained sample to all possible subsamples and construction of
regression estimates for each of these subsamples. The
decomposition for which the minimum of the general sum of
regression residuals is attained, is assumed to be the estimate of
a true decomposition of the whole samples of obtained observations
into subsamples with different regression regimes.

These methods turn out to be rather time consuming and have a low
power. For example, if there are only two regression regimes in an
obtained sample but we do not know this fact and are obliged to
try all possible subsamples up to the order 20, then many false
structural changes will be obtained.

In this paper we propose a new method of detection and estimation
of multiple change-points which is not based upon LSE of
regression parameters and computation of corresponding residuals.
This method is more effective and robust to possible inaccuracies
in specification of regression models.

Let us explain the idea of this method by the following example of
a multiple regression model (1) with deterministic predictors and
the row-matrix $\Pi(\vartheta,n)$. In other words, let
$\vartheta=(\theta_1,\theta_2,\dots,\theta_k),\,k\ge 1$ is an
unknown vector of change-point parameters such that $0\equiv
\theta_0<\beta\le \theta_1< \dots < \theta_k \le\alpha<
\theta_{k+1}\equiv 1$, where, as before, $\beta,\,\,\alpha$ are
known numbers, and the observations has the form
$$
y_n=\Pi^*(\vartheta,n)F(n/N)+\nu_n. \eqno(11)
$$
Here
$$
\Pi(\vartheta,n)=\suml_{i=1}^{k+1}\,a_i\,\mathbb{I}([\theta_{i-1}N]<n\le
[\theta_i N]),
$$
where $a_i\not=a_{i+1},i=1,2,\dots,k$ are unknown \emph{vectors},
$F(t)$ is a given vector-function (all assumptions and notations
see in Subsection 4.1.1).

Consider our main statistic (7). The mathematical expectation of
this statistic converges as $N\to\infty$ to the function
$$
m(t)=\intl_0^t\,F(s)F^*(s)\Pi(\vartheta,s)ds-A(t)I^{-1}
\intl_0^1\,F(s)F^*(s)\Pi(\vartheta,s)ds.
$$

In the situation when there is no change-points, i.e.,  the vector
of regression coefficients is constant on $[0,1]$, the vector
function $m(t)$ equals to zero for each $t\in [0,1]$. This
property of $m(t)$ makes it possible to effectively reject the
null hypothesis about the absence of change-points when they are
really present in an obtained sample.

Consider the following method of detection and estimation of
multiple change-points. Fix a small parameter
$\epsilon,\,\,\min(\beta,1-\alpha)>\epsilon>0$. The proposed
method consists of the following steps:

1. Compute statistic (7) by the data in the diapason of arguments $\mathcal{N}\df\left([\beta N],\dots,[\alpha N]\right)$. If
$\maxl_{n\in\mathcal{N}} \|\mathcal{Z}_{\scriptscriptstyle N}(n)\|^2> C$, where $C=C(N)$ is the decision threshold, then compute $nmax=\mathrm
{argmax}\|\mathcal{Z}_{\scriptscriptstyle N}(n)\|^2$, otherwise the sample is assumed to be stationary (without change-points).

2. Put $N^{'}=nmax-[\epsilon N]$ and compute statistic (7) by the
data in the diapason of arguments $\mathcal{N}^{'}\df\left([\beta
N],\dots,N^{'}\right)$ according to step 1. This cycle is repeated
until:

1) we obtain a stationary sub-sample in the diapason of data with
arguments $\left([\beta N],\dots,N^{'}\right)$, i.e.
$\maxl_{n\in\mathcal{N}^{'}} \|Z_{\scriptscriptstyle N^{'}}(n)\|^2
\le C(N^{'})$. Then we put $n(1)=N^{'}+[\epsilon N]$ as the
estimate of the first change-point and go to step 3.

or

2) we obtain a sample of the size $N^{'}\le [2\epsilon N]$. Then
we put $n(1)=N^{'}+[\epsilon N]$ as the estimate of the first
change-point and go to step 3.

3. Put $n^{'}=n(1)+[\epsilon N]$ and compute statistic (7) by the
data in the diapason of arguments $\left(n^{'},\dots,[\alpha
N]\right)$ (i.e. with the relative arguments $[1,\dots,[\alpha
N]-n^{'}+1]$) and do according to steps 1 and 2. The cycle is
repeated until we obtain a stationary sub-sample in the diapason
of data with arguments $[n^{'},\dots,nmax]$ or $nmax-n^{'}\le
[2\epsilon N]$. Then we put $n(2)=nmax$ as the estimate of the
next change-point. If $N-n(2) < [2\epsilon N]$ then stop,
otherwise repeat step 3 by the data in the diapason of arguments
$\left(n(2),\dots,[\alpha N]\right)$.

In this way we continue to compute the estimates $n(3),\dots$ of
change-points. As a result we obtain the series of estimates
$n(1),n(2),\dots$ of the true change-points $[\theta_1
N],\dots,[\theta_k N]$. The number $\hat k_N$ of these estimates
is determined by the quantity of stationary sub-samples
$$[1,\dots,n(1)],\dots,[n(i),\dots,n{(i+1)}],\dots,[n{(\hat
k_N)},\dots,N]$$.

The proposed method is based upon reduction to the case of only
one change-point and the properties of the matrix $m(t)$. The
crucial point of this method is the choice of the decision
threshold $C(N)$ which depends on the sample size $N$. Below we
give an explicit formula for computation of $C(N)$.

Let $\hat k_N$ be the estimate of the number of change-points in
the sample and $\hat{\vartheta}_N=(\theta_{N1},\dots,\theta_{N
{\hat k_{N}}})^*$ be the vector of estimated coordinates of
change-point parameters. The following theorem holds for model
(11).
\begin{thm}
Suppose assumptions of Theorem 3 are satisfied. Moreover, assume
that there exist $h>0,\,B>0$ such that for all $i=2,\dots,k+1$:
$$\begin{array}{ll}
& 0< \|A(\theta_{i-1},\theta_i)A^{-1}(\theta_{i-2},\theta_{i-1})\|\le h \\
& \|A(\theta_{i-1},\theta_i)(a_i-a_{i-1})\| \ge B >0,
\end{array}
$$
Then for sufficiently small $\delta>0$:
$$
\mathbf{P}\{(\hat k_N \ne k)\cup \{(\hat k_N = k)\cap (\maxl_{1\le
i\le k}\,|\hat\theta_{Ni}-\theta_i|> \delta )\} \} \le
C(\delta)\exp(-D(\delta)N),
$$
where constants $C(\delta)>0,D(\delta)>0$ do not depend on $N$.
\end{thm}
Analogous theorem can be proved also for stochastic predictors.

\bigskip
From theorem 5 it follows that the estimated number of
change-points converges almost surely to its unknown true value,
as well as estimated coordinates of unknown change-points converge
exponentially to their true values as the sample size tends to
infinity. Moreover, comparing results of theorem 2 and theorem 5
we conclude that the proposed method of detection and estimation
of multiple change-points is asymptotically optimal by the order
of convergence of estimated change-point parameters to their true
values.

The proof of theorem 5 is given in the Appendix E.

\subsection{A variant of the limit distribution theorem for the decision statistic under the null hypothesis}

For practical applications of the proposed method and, in
particular, for the rational choice of the decision threshold
$C(N)$, we need to study the limit distribution of the decision
statistic under the null hypothesis.

Let us formulate a variant of the limit theorem for the simple
case of unique change-point, deterministic predictors,
statistically independent noises $\nu_n$, and the one-dimensional
dependent variable $y_n$.

Suppose there exists a continuous function $g(t),\,0\le t\le 1$
such that $\mathbf{E}_{\theta}\,\nu_n^2=g^2(n/N)$.

Put
$$
\begin{array}{ll}
&\sigma_i^2=\dpfrac 1t\,\intl_0^t\,f_i^2(s)g^2(s)ds,\,i=1,\dots,K\\[2mm]
&G(t)=(\sigma_1(t),\dots,\sigma_{\scriptscriptstyle
K}(t))^*,\,\,\mathbf{Z}(t)=G(t)W(t),\,\,U(t)=\mathbf{Z}(t)-A(t)I^{-1}\mathbf{Z}(1),
\end{array}
$$
where $W(t)$ is the standard Wiener process, $A(t),\, I$ are the
above defined matrices (see Subsection 4.1.1).

Consider our main statistic, the vector process
$\mathcal{Z}_{\scriptscriptstyle
N}(t)=\mathcal{Z}_{\scriptscriptstyle N}([Nt])$ (see (7)). Then
for any $\theta \in [\beta,\alpha]$, the vector process $\sqrt{N}
(\mathcal{Z}_{\scriptscriptstyle
N}(t)-\mathbf{E}_{\theta}\,\mathcal{Z}_{\scriptscriptstyle N}(t))$
weakly converges to the vector process $U(t)$ in the Skorokhod
space $D^K\,[\beta,\alpha]$ (see Brodsky, Darkhovsky (2000)). In
particular, under the null hypothesis, the weak convergence is
valid at $[0,1]$.

Therefore, we have the following
\begin{thm}
$$
\liml_{N\to\infty}\,\mathbf{P}_0\{\sqrt{N}\maxl_{t\in
[0,1]}\,\|\mathcal{Z}_{\scriptscriptstyle N}(t)\| >C\}=
\mathbf{P}_0 \{\maxl_{t\in [0,1]}\,\|U(t)\|> C\} \eqno(12)
$$
(here we use the Euclidean norm for vectors).
\end{thm}
The vector $U(t)$ is Gaussian with zero mean and the following
$K\times K$ correlation matrix $D(t)$:
$$
D(t)=t\left[G(t)G^*(t)-G(t)G^*(1)I^{-1}A(t)-A(t)G(1)G^*(t)\right]
+A(t)I^{-1}G(1)G^*(1)I^{-1}A(t).
$$
Therefore, we have the following equality by distribution
$$
U(t)=\sqrt{D(t)}\zeta \eqno(13)
$$
where $\zeta=(\zeta_1,\dots,\zeta_{\scriptscriptstyle K})^*$ is
the standard Gaussian vector.

Taking (13) into account, we get
$$
\maxl_{0\le t\le 1}\|U(t)\|=\maxl_{0\le t\le 1}\sqrt{\suml_{i=1}^K
d_i^2(t)\zeta_i^2}\df\rho(\zeta), \eqno(14)
$$
where $d_i^2(t)$ are eigenvalues of the matrix $D(t)$. The
function $\rho(\zeta)$ can be explicitly calculated for any given
family of functions $F(t), g(t)$.

Therefore, from (14) we have
$$
\mathbf{P}_0\{\maxl_{0\le t\le
1}\|U(t)\|>C\}=\intl_{\{u:\rho(u)>C\}} \varphi(u)du, \eqno(15)
$$
where $\varphi(u)$ is the density of the standard Gaussian
distribution.

From (12) and (15) we can conclude that type 1 error goes to zero
as $\exp(-const \,NC^2)$ for the proposed method. This fact allows
us to choose the decision threshold. Note that the same
asymptotical order can be obtained from corollary 2 (see
Subsection 4.1.1). For independent noises we have
$$
\mathbf{P}_0 \{\maxl_{[\beta N]\le n\le N}
 \|{\cal Z}_{\scriptscriptstyle N}(n)\|^2>C \}\le
\left \{
\begin{array}{ll}
& \exp \left(-\dpfrac {TN\mathbb{C}\beta}{4R}\right), \ \
\mathbb{C} >gT \\
& \exp \left(-\dpfrac {N\beta \mathbb{C}^2}{4R^2g
m_0(\mathbb{C})}\right), \ \mathbb{C} \le gT,
\end{array}
\right.
$$
(the notations see in Subsection 4.1.1).

Therefore, we conclude that type 1 error
$\alpha_{\scriptscriptstyle N}$ goes to zero exponentially as
$N\to\infty$ for the proposed method.

So, the threshold can be calculated from the relation
$$
C=C(N)=\dpfrac 1{\sqrt{N}}\,|\ln\alpha_{\scriptscriptstyle
N}|\,\lambda,
$$
where $\lambda$ is a certain calibration parameter which depends
on variations of predictors, dispersions of noises and
characteristics of their statistical dependence.

A more close study allows us to obtain the following practical
formula for the decision threshold $C=C(N)$:
$$
C(N)=\dpfrac {\Big({\maxl_i\sigma_i^2}\,\cdot{\maxl_i \maxl_{0\le
t\le 1}f_i^2(t)}\Big)^{1/2}}{\sqrt{N}}\,\lambda,
$$
where $\sigma_i^2$ is the dispersion of $\nu_i$ and $\lambda>0$ is
the calibration parameter.

\section{Experiments}

In this section we present results of a simulation study of the
proposed method in comparison with other well known tests. The
following methods are most often used for detection of structural
changes in regression models:

- The Chow test most often used in econometric packages;

- The CUSUM (cumulative sums) test based upon recursive regression
residuals (Brown, Durbin, Evans, 1975);

- The CUSUM test based upon residuals of ordinary least squares
method (OLS CUSUM test, Ploberger, Kramer, 1992);

- Fluctuation test (Ploberger, Kramer, Kontrus, 1989)

- Wald test (Andrews, 1993, Andrews, Ploberger, 1994)

- LM тест (Lagrange Multilpier test, Andrews, 1993).

However, it is well known (see, e.g., Maddala and Kim (1998)) that
the Wald test (together with  the QMLE - quasi-maximum likelihood
estimation test) is the best and most often used for detection of
changes in regression models because it has the best
characteristics of power and accuracy of change-point estimation.

The Wald test statistic is defined as follows:
$$
Sup W = \maxl_{1\le m\le N}\,N[\dpfrac
{S(N)-S_1(m)-S_2(N-m)}{S_1(m)+S_2(N-m)}],
$$
where $S(N)$ is the sum of regression residuals constructed by the
whole sample of the size $N$; $S_1(m)$ is the sum of regression
residuals constructed by the sub-sample of the first $m$
observations; $S_2(N-m)$ is the sum of residuals of the regression
model constructed by the last $N-m$ observations.

It is natural to define the estimate of the change point as
$n_0\in arg\sup\,W$, and the corresponding estimate of the
change-point parameter $\hat\theta_N=n_0/N$.

Comparison of characteristics of different methods is carried out
in the following way. First, methods are 'equalized' by the value
of type 1 error by means of choice of the corresponding decision
thresholds. In practice, for this purpose we use experiments with
stationary samples (without structural changes) in which the
95-percent quantiles of the variation series of the decision
statistics are computed (see below, table 1). Second, for the
chosen sample sizes and decision thresholds, experiments with
non-stationary samples are performed in which we compute estimates
of the type 2 error probability and instants of change-points (see
tables 2 and 4). The method of change-point detection  'a' is
preferable w.r.t. the method "b" if for the same values of the
type 1 error, it gives lower estimates of the type 2 error and the
error of change-point estimation.

\subsection{Deterministic regression plan}
We compared characteristics of our method with those of the Wald
test using the following regression model with deterministic
predictors:
$$
y_i=c_0+c_1\,x_i+\xi_i,   \quad i=1,\dots, N       \eqno (16)
$$
where $(x_1,\dots,x_N)^*$ is the vector of deterministic
predictors; $\{\xi_i\}$ is the Gaussian noise sequence with zero
mean and unit variance; $c_0,c_1$ are regresson coefficients which
change at the instant $n_0=[\theta\,N],\,0< \theta<1$.

The number of independent trials of each experiment was equal to
k=2000. The estimates of decision thresholds were obtained as
follows. For each stationary sample, the 95-percent and 99-percent
quantiles of the variation series of maximums of the decision
statistic were computed in 2000 trials. These quantiles were then
assumed to be estimates of the decision thresholds for 5-percent
and 1-percent error level, respectively.

The values of the threshold $C$ given in table 1, were used as
decision bounds for the confidence probability 95 percent in
experiments with non-stationary regression models. The following
cases were considered:

- before the change-point: $c_0=0,\,c_1=1$

- after the change-point: $c_0=\delta,\,c_1=1$.

In experiments the parameter $\delta$ and the sample size $N$ were
changed. The following characteristics of the proposed method were
estimated:

- The empirical estimate of decision threshold $C$ (more exactly,
the empirical estimate of
$\maxl_n\|\mathcal{Z}_{\scriptscriptstyle N}(n)\|$);

- The empirical estimate of type 2 error probability $\hat w_N$;

- The empirical estimate of the change-point parameter
$\hat\theta_N$.

Results obtained for the Wald test are given in the following
tables.

{\bf Table 1. Estimation of the decision thresholds for the Wald
test for different sample sizes}

\medskip

\begin{tabular}{|c|c|c|c|c|c|c|c|c|}
\hline
$N$ & 100 & 200 & 300 & 400 & 500 & 700 & 1000 & 1200 \\
\hline
$p=0.95$ & 10.10 & 8.09 & 9.59 & 8.66 & 8.12 & 7.62 & 7.51 & 7.43\\
\hline
$p=0.99$ & 12.60 & 10.88 & 14.14 & 12.10 & 12.20 & 9.97 & 11.68 & 10.02\\
\hline
\end{tabular}

\bigskip
{\bf Table 2. Estimation of the change-point parameter
$\theta=0.30$ by the Wald test}

\medskip

\begin{tabular}{|c|c|c|c|c|c|c|}
\hline
\multicolumn{2}{|c|}{$N$} & 300 & 400 & 500 & 700 & 1000 \\
\hline
$\delta=0.3$ & $C$ & 5.63 & 6.76 & 8.24 & 9.77 & 12.09 \\
\cline{2-7}
  & $\hat w_N$ & 0.83 & 0.71 & 0.59 & 0.46 & 0.32 \\
\cline{2-7}
 & $\hat\theta_N$ & 0.29 & 0.25 & 0.22 & 0.19 & 0.20 \\
\hline
$\delta=0.4$ & $C$ & 9.65 & 10.20 & 11.88 & 15.27 & 19.32 \\
\cline{2-7}
  & $\hat w_N$  & 0.56  & 0.47 & 0.34 & 0.23 & 0.18 \\
\cline{2-7}
 & $\hat\theta_N$ & 0.28 & 0.25 & 0.22 & 0.20 & 0.23 \\
\hline
\end{tabular}

\bigskip
The same model was studied with the help of the method proposed in
this paper.

1) Decision thresholds

In the first series of experiments, model (16) with constant
coefficients $c_0=0,\,c_1=1$ was used. The following results were
obtained.

{\bf Table 3. Estimation of the decision thresholds}

\medskip

\begin{tabular}{|c|c|c|c|c|c|c|c|c|}
\hline
$N$ & 100 & 200 & 300 & 400 & 500 & 700 & 1000 & 1200 \\
\hline
$p=0.95$ & 0.401 & 0.257 & 0.202 & 0.182 & 0.150 & 0.125 & 0.103 & 0.081\\
\hline
$p=0.99$ & 0.450 & 0.300 & 0.247 & 0.211 & 0.187 & 0.162 & 0.138 & 0.102\\
\hline
\end{tabular}

\bigskip
2) The estimates of the change-point parameter

{\bf Table 4. Results of estimation of the change-point parameter
$\theta=0.30$}

\medskip

\begin{tabular}{|c|c|c|c|c|c|c|}
\hline
\multicolumn{2}{|c|}{$N$} & 300 & 400 & 500 & 700 & 1000 \\
\hline
$\delta=0.3$ & $C$ & 0.179 & 0.177 & 0.168 & 0.157 & 0.151 \\
\cline{2-7}
  & $\hat w_N$ & 0.64 & 0.55 & 0.33 & 0.13 & 0.03 \\
\cline{2-7}
 & $\hat\theta_N$ & 0.340 & 0.322 & 0.332 & 0.324 & 0.307 \\
\hline
$\delta=0.4$ & $C$ & 0.220 & 0.211 & 0.208 & 0.195 & 0.192 \\
\cline{2-7}
  & $\hat w_N$  & 0.28  & 0.24 & 0.11 & 0.02 & 0.005 \\
\cline{2-7}
 & $\hat\theta_N$ & 0.315 & 0.312 & 0.308 & 0.305 & 0.304 \\
\hline
\end{tabular}

\bigskip
{\bf Table 5. Results of estimation of the change-point parameter
$\theta=0.50$}

\medskip

\begin{tabular}{|c|c|c|c|c|c|c|}
\hline
\multicolumn{2}{|c|}{$N$} & 300 & 400 & 500 & 700 & 1000 \\
\hline
$\delta=0.3$ & $C$ & 0.194 & 0.184 & 0.175 & 0.168 & 0.164 \\
\cline{2-7}
  & $\hat w_N$ & 0.62 & 0.50 & 0.25 & 0.05 & 0.01 \\
\cline{2-7}
 & $\hat\theta_N$ & 0.456 & 0.485 & 0.501 & 0.502 & 0.499 \\
\hline
$\delta=0.4$ & $C$ & 0.231 & 0.221 & 0.215 & 0.214 & 0.211 \\
\cline{2-7}
  & $\hat w_N$  & 0.26  & 0.22 & 0.003 & 0.02 & 0 \\
\cline{2-7}
& $\hat\theta_N$ & 0.495 & 0.495 & 0.489 & 0.501 & 0.499 \\
\hline
\end{tabular}

\bigskip
Comparing results from tables 2 and 4, we conclude that type 2
error estimates for our method are  lower than for the Wald test,
and the error of estimation for our method is much lower than for
the Wald test. Therefore, we conclude that our method is
essentially better by the main performance characteristics of
change-point detection than the Wald test, and so, we conclude
that the proposed method is one of the most effective among all
known tests for detection and estimation of structural changes in
regression models.

Comparing results from table 4 and 5, we can conclude that the
quality of estimation of the change-point parameter $\theta$
depends on its location on the segment $[0,1]$: estimation of
$\theta$ which is closer to the bounds of the segment $[0,1]$ is
more difficult.

\bigskip
In next two subsections we investigate our methods.
\subsection{Stochastic regression plan}
In this series of experiments the following model of observations
was used:
$$
y_i=c_0+c_1\,x_i+\xi_i,   \quad i=1,\dots, N
$$
where $(x_1,\dots,x_N)^{*}$ is a stationary random sequence of the
following type:
$$
x_i=\rho x_{i-1}+\eta_i, \quad i=1,\dots,N,\; x_0\equiv 0,
$$
$\{\xi_i,\,\eta_i\}$ is the sequence of independent Gaussian
r.v.'s with zero mean and unit dispersion; $c_0,c_1$ are
regression coefficients which change at the instant
$n_0=[\theta\,N],\,0< \theta<1$; $|\rho| <1$.

1) Estimation of decision thresholds

In the first series of tests decision thresholds were estimated.
For this purpose, stationary sequences (without change-points)
were used: $c_0=0,\,c_1=1,\rho=0.3$. The following results were
obtained.

\bigskip
{\bf Table 6. Estimation of decision thresholds (the case of
stochastic predictors)}

\begin{tabular}{|c|c|c|c|c|c|c|c|c|}
\hline
$N$ & 100 & 200 & 300 & 400 & 500 & 700 & 1000 & 1200 \\
\hline
$p=0.95$ & 0.355 & 0.291 & 0.230 & 0.188 & 0.150 & 0.132 & 0.103 & 0.082\\
\hline
$p=0.99$ & 0.401 & 0.332 & 0.273 & 0.218 & 0.192 & 0.171 & 0.141 & 0.100\\
\hline
\end{tabular}

\bigskip
2) Estimation of the change-point parameter

In the following series of experiments a model with a structural
change in the regression coefficients was used:

- before the change-point: $c_0=0,\,c_1=1$

- after the change-point: $c_0=0,\,c_1=1.3$.

Results obtained are presented in table 7.

\bigskip
{\bf Table 7. Estimation of change-point parameters (the case of
stochastic predictors)}

\begin{tabular}{|c|c|c|c|c|c|}
\hline
\multicolumn{2}{|c|}{$N$} & 500 & 700 & 1000 & 1200 \\
\hline
$\theta=0.5$ & $C$ & 0.167 & 0.157 & 0.152 & 0.152 \\
\cline{2-6}
  & $\hat w_N$ & 0.32 & 0.21 & 0.02 & 0 \\
\cline{2-6}
 & $\hat\theta_N$ & 0.481 & 0.495 & 0.498 & 0.499 \\
\hline
$\theta=0.3$ & $C$ & 0.156 & 0.148 & 0.142 & 0.140 \\
\cline{2-6}
  & $\hat w_N$  & 0.45  & 0.30 & 0.03 & 0 \\
\cline{2-6}
& $\hat\theta_N$ & 0.312 & 0.310 & 0.308 & 0.301 \\
\hline
\end{tabular}

\subsection{Multiple structural changes in multivariate systems}
The following multivariate system was used:
$$\begin{array}{ll}
& y_i=c_0+c_1 y_{i-1}+c_2 z_{i-1}+c_3 x_i+\epsilon_i \\
& z_i=d_0+d_1 y_i+d_2 x_i+\xi_i \\
& x_i=0.5 x_{i-1}+\nu_i \\
& \epsilon_i=0.3\epsilon_{i-1}+\eta_i,
\end{array}
$$
where $\xi_i,\nu_i,\eta_i,\;i=1,2,\dots$ are independent standard
Gaussian random variables.

Here $(y_i,z_i)^*$ is the vector of endogenous variables, $x_i$ is
the vector of exogenous variables, $(y_{i-1},z_{i-1},x_i)^*$ - the
vector of pre-determined variables of the considered system.

Dynamics of this system is characterized by the following vector
of coefficients: $\mathbf {u}=[c_0\; c_1\; c_2\; c_3\; d_0\; d_1\;
d_2]$. The initial vector of coefficients is $[0.1\; 0.5\; 0.3\;
0.7\; 0.2\; 0.4\; 0.6]$. The first structural change occurs at the
instant $\theta_1=0.3$. The vector of coefficients $\mathbf {u}$
changes into $[0.1\; 0.5\; 0\; 0.7\; 0.2\; 0.4\; 0.6]$. The second
structural change occurs at the instant $\theta_2=0.7$. Then the
vector $\mathbf {u}$ changes into $[0.1\; 0.5\; 0\; 0.7\; 0.2\;
0.4\; 0.9]$.

In the first series of tests the decision threshold $C$ was
estimated. For this purpose, the model with the initial vector of
coefficients $\mathbf {u}$ and without change-points was used. In
2000 independent trials the maximums of the decision statistic
were computed and the variation series of these maximum was
constructed. Then the 95-percent and  the 99-percent quantiles of
this series were computed. These values are presented in table 8.

{\bf Table 8. Estimation of decision thresholds (the case of a
multivariate system)}

\begin{tabular}{|c|c|c|c|c|c|c|c|c|}
\hline
$N$ & 200 & 400 & 500 & 700 & 900 & 1000 & 1200 & 1500 \\
\hline
$p=0.95$ & 0.28 & 0.20 & 0.19 & 0.18 & 0.16 & 0.15 & 0.145 & 0.14\\
\hline
$p=0.99$ & 0.36 & 0.33 & 0.28 & 0.24 & 0.23 & 0.21 & 0.19 & 0.17\\
\hline
\end{tabular}

\bigskip
The computed 95-percent quantiles were assumed to be the decision
thresholds for the corresponding sample volumes.

In the next series of tests non-stationary samples with multiple
change-points were used. The true number of change-points was
equal to $p=2$, the coordinates of these change-points were
$\theta_1=0.3$ and $\theta_2=0.7$. In table 9 the following
performance characteristics are given:

- $w$ is the estimate of the probability
$\mathbf{P}_{\theta}\{\hat p_{\scriptscriptstyle N}\ne p\}$ in
2000 independent trials, where $\hat{p}_{scriptscriptstyle N}$ is
the estimate of the number of change-points in the data.

- $\Delta$ is the estimation error on condition that $\hat p_N=
p$, i.e.  $\Delta=\sqrt
{\sum_{i=1}^p\,(\hat\theta_i-\theta_i)^2}$.

{\bf Table 9. Estimation of change-point parameters (the case of a
multivariate system)}

\begin{tabular}{|c|c|c|c|c|c|c|c|c|}
\hline
$N$ & 200 & 400 & 500 & 700 & 900 & 1000 & 1200 & 1500 \\
\hline
$ w $ & 0.96 & 0.54 & 0.39 & 0.21 & 0.04 & 0.03 & 0.02 & 0.01\\
\hline
$ \Delta $ & 0.02 & 0.05 & 0.04 & 0.02 & 0.03 & 0.02 & 0.01 & 0.005 \\
\hline
\end{tabular}

\section{Conclusions}

In this paper the following main results were obtained:

1. The general statement of the retrospective change-point
detection and estimation problem in multivariate linear systems is
given (both one change-point and multiple change-point problems,
both independent and dependent sequences of observations)

2. The prior lower bounds are proved for the main performance
characteristic in retrospective change-point detection and
estimation: {\it the probability of the error of change-point
estimation}, in different contexts of change-point estimation:
from one change-point in multi-factor linear regressions with
deterministic and stochastic regression plans, to multiple
change-point problems in multivariate linear models.

3. A new method is proposed for the problem of retrospective
change-point detection and estimation in multivariate linear
systems. The main performance characteristics of this method: type
1 and type 2 errors, the error of change-point estimation, are
studied theoretically. We prove that the proposed method is {\it
asymptotically optimal} by the order of convergence of the
change-point estimate to its true value as the sample size tends
to infinity.

4. For the problem of multiple change-point detection and
estimation, we propose a general setup in which both {\it the
number of change-points and their coordinates in the sample are
unknown}. For this problem statement, a new method is proposed
which gives consistent estimates of an unknown number of
change-points and their coordinates. This method is also
asymptotically optimal by the order of convergence of these
estimates to true change-point parameters.

5. A simulation study of characteristics of the proposed method
for finite sample sizes is performed. The main goals of this study
are as follows: to compare performance characteristics of the
proposed method with characteristics of other well known methods
of change-point detection in linear regression models: the Wald
test, the Chow test, the CUSUM tests with ordinary and recursive
regression residuals, the fluctuation test; to consider more
general multivariate linear models and performance characteristics
of the proposed method in these multivariate models. The main
conclusion: performance characteristics of the proposed method are
no worse but often even better than those of well known
change-point tests.

\bigskip
{\Large \bf Appendix. Proofs of theorems}

\appendix
\section{Proof of Theorem 1}

Using notations (3)-(4), put
$$
\mathcal{M}(\Delta)=\{T(\Delta):T(\Delta)=\{T_{\scriptscriptstyle
N}(\Delta)\}_{\scriptscriptstyle N=1}^{\infty}\}
$$
This is the set of all sequences of the elements
$T_{\scriptscriptstyle
N}(\Delta)\in\mathcal{M}_{\scriptscriptstyle N}(\Delta)$. Consider
also the collection of all consistent estimates of the parameter
$\theta\in\Delta$, i.e.,
$$
\tilde{\mathcal{M}}(\Delta)=\{T(\Delta)\in\mathcal{M}(\Delta):\liml_{N\to\infty}\mathbf{P}_{\theta}(|T_{\scriptscriptstyle
N}(\Delta)-
\theta|>\epsilon)=0,\,\forall\theta\in\Delta,\,\forall\epsilon>0\}
$$
Under the assumption of Theorem 1, the set
$\tilde{\mathcal{M}}([a,b])$ is \emph{non-empty} for any
$0<a<b<1$. Indeed, consider the sequence
$y_n=\ln\dpfrac{f_0(x_n,n/N)}{f_1(x_n,n/N)}$. Due to the
assumption, $\mathbf{E}_{\theta}y_n\ge\delta>0$ before the
change-point $\theta,\,a\le\theta\le b$, and less than $(-\delta)$
after the change-point. Now, using the same idea as in Brodsky and
Darkhovsky (2000), it is easy to construct the consistent estimate
of the change-point.

Further, without loss of generality we can consider only
consistent estimates of the change-point parameter $\theta$,
because for non-consistent estimates the probability of the error
of estimation does not converge to zero and the considered
inequality is satisfied trivially.

Let $\hat\theta_N$ be some consistent estimate of the change-point
parameter $\theta$ constructed by the sample
$X^N=\{x_1,\dots,x_{\scriptscriptstyle N}\}$. Consider the random
variable $\lambda_{\scriptscriptstyle
N}=\lambda_{\scriptscriptstyle N}(x_1,\dots,x_{\scriptscriptstyle
N})= \mathbb{I}\{|\hat\theta_N-\theta|> \epsilon \}$.

Under the change-point parameter $\theta$, the likelihood function
for the sample $X^N$ can be written as follows:
$$
f(X^N,\theta)=\prodl_{i=1}^{[\theta N]}\,f_0(x_i,i/N)\cdot
\prodl_{i=[\theta N]+1}^N\, f_1(x_i,i/N).
$$

We have for any $d>0$ and $0< \epsilon < \epsilon^{'}$:
$$
\begin{array}{ll}
&\mathbf{P}_{\theta}\{|\hat\theta_N-\theta|> \epsilon
\}=\mathbf{E}_{\theta}\lambda_{\scriptscriptstyle N} \ge
\mathbf{E}_{\theta}(\lambda
\mathbb{I}(f(X^N,\theta+\epsilon^{'})/f(X^N,\theta)
< e^d))\ge\\[2mm]
&\ge e^{-d} \left(\mathbf{E}_{\theta +\epsilon\prime}(\lambda
_{\scriptscriptstyle N} \mathbb{I}(f(X^{N},\theta
+\epsilon\prime)/f(X^{N},\theta )<
e^{d}\}\right)\ge\\[2mm]
& e^{-d} \left(\mathbf{P}_{\theta +\epsilon\prime}\{| \theta _{N}
- \theta| >\epsilon \} - \mathbf{P}_{\theta
+\epsilon\prime}\{f(X^{N},\theta +\epsilon\prime)/f(X^{N},\theta )
\ge e^{d}\}\right)
\end{array}
$$
(here we used the elementary inequality $\mathbf{P}(AB) \ge
\mathbf{P}(A) - \mathbf{P}(\Omega \backslash B)$).

Consider the probabilities in the right-hand side of the last
inequality. Since $\theta _{N}$ is a consistent estimate of
$\theta$, we have $\mathbf{P}_{\theta +\epsilon\prime}\{| \theta
_{N} - \theta| > \epsilon \}\rightarrow 1$ as $N\rightarrow\infty
$. For estimation of the second probability, we take into account
that
$$
\ln \left(f(X^{N},\theta +\epsilon\prime)/f(X^{N},\theta )\right)
= \suml^{[(\theta +\epsilon\prime)N]}_{i=[\theta N]+1} \ln
\left(f_{0}(x_i,i/N)/f_{1}(x_i,i/N)\right)
$$
Therefore,
$$
\begin{array}{ll}
&\mathbf{E}_{\theta +\epsilon\prime}\ln \left(f(X^{N},\theta
+\epsilon\prime)/f(X^{N},\theta )\right)=\\[2mm]
&= N\,\intl_{\theta}^{\theta+\epsilon{'}}\,\mathbf{E}_0\ln\dpfrac
{f_0(x,t)}{f_1(x,t)}dt+O(1).
\end{array}
$$

Then
$$
\begin{array}{ll}
&\mathbf{P}_{\theta +\epsilon{'}}\{f(X^{N},\theta
+\epsilon{'})/f(X^{N},\theta ) \ge  e^{d}\} =\\[2mm]
&=\mathbf{P}_{\theta +\epsilon{'}}\left\{\suml_{i=[\theta
N]+1}^{[(\theta+\epsilon^{'})N]}\,
(\ln(f_0(x_i,i/N)/f_1(x_i,i/N))-\mathbf{E}_0\ln
\left(f_0(x_i,i/N)/f_1(x_i,i/N)\right)\right.\\[2mm]
&\left.\ge
d-N\intl_{\theta}^{\theta+\epsilon{'}}\,\mathbf{E}_0\ln\dpfrac
{f_0(x,t)}{f_1(x,t)}dt+O(1)\right\}
\end{array}
$$
Put
 $d=d_1(N)=N(\intl_{\theta}^{\theta+\epsilon{'}}\,\mathbf{E}_0\ln\dpfrac
{f_0(x,t)}{f_1(x,t)}dt+ \delta)$ for some $\delta>0$ and use the
law of large numbers which holds due to existence of
$\mathbf{E}_0\ln\dpfrac{f_0(x,t)}{f_1(x,t)}$. Then we obtain
$$
\mathbf{P}_{\theta +\epsilon{'}}\{f(X^{N},\theta
+\epsilon{'})/f(X^{N},\theta ) \ge  e^{d_1(N)}\}\to 0
$$
as $N\to\infty$.

The same considerations for $d=d_{2}(N)=
N(\intl_{\theta-\epsilon{'}}^{\theta}\,\mathbf{E}_1\ln\dpfrac
{f_1(x,t)}{f_0(x,t)}dt+\delta)$ yield
$$
\mathbf{P}_{\theta-\epsilon{'}}\left\{f(X^{N},\theta
-\epsilon{'})/f(X^{N},\theta ) \ge  e^{d_{2}(N)}\right\}
\rightarrow 0
$$
as $N\rightarrow\infty $.

Therefore,
$$
\mathbf{P}_{\theta}\{|\hat\theta_N-\theta|> \epsilon\}\ge
(1-o(1))\max (e^{-d_1(N)},e^{-d_2(N)}).
$$

It follows from here
$$
\liminf_{N\to\infty}\,N^{-1}\ln\infl_{\hat\theta_N\in
{\mathcal{M}}_N}\, \mathbf{P}_{\theta}\{|\hat\theta_N-\theta|>
\epsilon \} \ge -\min\left
(\intl_{\theta}^{\theta+\epsilon^{'}}\,J_0(t)dt,\,\intl_{\theta-
\epsilon^{'}}^{\theta}\,J_1(t)dt\right)-\delta.
$$

Note that the left-hand side of this inequality does not depend on
the parameters $\delta,\,\epsilon^{'}$, and the right-hand side
exists for each $\delta>0,\,\theta\wedge
(1-\theta)>\epsilon{'}>\epsilon
>0$. From the continuity assumption for the functions $J_0(\cdot),
J_1(\cdot)$, we conclude that our result follows after taking the
limits of both sides of this inequality as $\delta \to 0$ and
$\epsilon^{'}\to \epsilon$.

\section{Proof of Theorem 2} We will use notations (3)-(4) and
(6). Let $x\in \mathbb{R}^{p}, y\in \mathbb{R}^{q}$, и $m=\max
(p,q)$. Define the following natural immersions:
$$
\mathrm{im}_{x}: \mathbb{R}^{p}\rightarrow\mathbb{R}^{m},\quad
\tilde {x}=\mathrm{im}_{x}\,x,\quad \mathrm{im}_{y}:
\mathbb{R}^{q}\rightarrow\mathbb{R}^{m},\quad \tilde {y}= \mathrm
{im}_{y}\,y
$$
(all lacking components are substituted by zeros) and put:
$$
\mathit{dist}(x,y) = \| \tilde{x}-\tilde{ y}\|^{(m)}
$$
(here we use the $\|\cdot\|_{\infty}$-norm for vector
$x=(x_1,\dots,x_p)$, i.e., $\|x\|^{(p)}=\maxl_{1\le i\le
p}|x_i|$).

Consider
$$
\begin{array}{ll}
&\liminf_{N\to\infty}N^{-1}\ln\infl_{\displaystyle\vartheta_{N}\in
\mathcal{M}_{N}
(\mathcal{D}^{\star})}\supl_{\displaystyle\vartheta\in
\mathcal{D}_{k}} \left\{\mathbf{P}_{\vartheta } (\vartheta _{N}\in
\mathcal{D}_{k},\| \vartheta _{N}-\vartheta \|^{(k)}>
\epsilon )\right.\\[3mm]
&\left.+ \mathbf{P}_{\vartheta } (\vartheta _{N}\not\in
\mathcal{D}_{k})\right \}
\end{array}
\eqno(B.1)
$$

Note that for $\epsilon  < \delta $, any estimate $\vartheta
_{N}\in \mathcal{M}_{N}(\mathcal{D}^{\star})$, and any $\vartheta
\in \mathcal{D}_{k}$, the following relationships between events
hold:
$$
\begin{array}{ll}
&\left(\mathit{dist}(\vartheta _{N},\vartheta)>\epsilon\right)
=\left(\vartheta_{N}\in \mathcal{D}_{k},\| \vartheta
_{N}-\vartheta \|^{(k)}
>\epsilon\right)\bigcup\left(\vartheta _{N}\not\in \mathcal{D}_{k},
\mathit{dist}(\vartheta _{N},\vartheta ) > \epsilon \right) =\\[3mm]
&= \left(\vartheta _{N}\in \mathfrak{D}_{k},\| \vartheta
_{N}-\vartheta \|^ {(k)}>\epsilon \right)\bigcup \left(\vartheta
_{N}\not\in \mathcal{D}_{k} \right).
\end{array}
$$

Here we used the fact that from the definition of \emph{ dist} and
the condition $(\vartheta _{N}\not\in \mathcal{D}_{k})$ it follows
that $(\mathit{dist}(\vartheta_{N},\vartheta ) > \delta )$, and
this condition yields $\mathit{dist}(\vartheta _{N},\vartheta ) >
\epsilon )$ for $\epsilon  < \delta $.

Thus, we need to estimate the probability $\mathbf{P}_{\vartheta
}\big( \mathit{dist}(\vartheta _{N},\vartheta ) > \epsilon \big)$.

First, note that the set $\tilde{\mathcal{M}}(\mathcal{D}_{k})$ of
all consistent estimates of the parameter $\vartheta  \in
\mathfrak{D}_{k}$ is \emph{non-empty}. This fact follows from
assumption ii) of the Theorem 2 and the same considerations as in
proof of Theorem 1.

Second, remark that the infimum in (B.1) can be taken only on the
set $\mathcal{M}_{N}(\mathcal{D}_{k})$. In fact, let $\vartheta
^{\star}_{N}\in \mathcal{M}_{N}(\mathcal{D}^{\star})$ belongs to
$\textrm{arg}\inf$ of the left-hand side of this inequality, i.e.,
$$
\begin{array}{ll}
&\infl_{\displaystyle\vartheta_{N}\in
\mathcal{M}_{N}(\mathcal{D}^{\star})}
\supl_{\displaystyle\vartheta\in
\mathfrak{D}_{k}}\mathbf{P}_{\vartheta } \{\mathit{dist}(\vartheta
_{N},\vartheta ) >\epsilon \}\\[3mm]
&=\supl_{\displaystyle\vartheta\in
\mathcal{D}_{k}}\mathbf{P}_{\vartheta }
\{\mathit{dist}(\vartheta^{\star}_{N},\vartheta ) >\epsilon \}
\end{array}
$$
(without loss of generality we suppose that the infimum is
attainable). Then consider the following element
$\hat{\vartheta}_{N}$ of the set
$\mathcal{M}_{N}(\mathcal{D}_{k})$:
$$
\hat{\vartheta}_{N}=
\vartheta^{\star}_{N}\mathbb{I}(\vartheta^{\star}_{N} \in
\mathcal{D}_{k}) + \Gamma _{N}\mathbb{I}(\vartheta^{\star}_{N}
\not\in \mathcal{D}_{k})
$$
where $\Gamma _{N}$ is the element of the set
$\mathcal{M}_{N}(\mathcal{D}_{k})$ such that
$$
\sup_{\displaystyle\vartheta\in
\mathcal{D}_{k}}\mathbf{P}_{\vartheta } \{\mathit{dist}(\Gamma
_{N},\vartheta ) < \epsilon /2\} \ge  1-\varkappa
$$
for some fixed $\varkappa >0$. Such elements exist in
$\mathcal{M}_{N}(\mathcal{D}_{k})$ (for large enough $N$), because
this set contains consistent estimates.

By definition, $\hat{\vartheta}_{N} \in
\mathcal{M}_{N}(\mathfrak{D}_{k})$ and for each $\vartheta  \in
\mathcal{D}_{k}$,
$$
\mathbf{P}_{\vartheta
}\{\mathit{dist}(\hat{\vartheta}_{N},\vartheta ) >\epsilon \}=
 \mathbf{P}_{\vartheta }\{\mathit{dist}(\vartheta^{\star}_{N},\vartheta ) >\epsilon \}
+ \mathbf{P}_{\vartheta }\{\mathit{dist}(\Gamma _{N},\vartheta )
>\epsilon \}.
$$
Therefore,
$$
\begin{array}{ll}
&\supl_{\displaystyle\vartheta\in
\mathcal{D}_{k}}\mathbf{P}_{\vartheta }
\{\mathit{dist}(\hat{\vartheta}_{N},\vartheta ) >\epsilon \} \le
\supl_{\displaystyle\vartheta\in \mathcal{D}_{k}}\mathbf{P}_{\vartheta } \{\mathit{dist}(\vartheta^{\star}_{N},\vartheta ) >\epsilon \}\\
&+\supl_{\displaystyle\vartheta\in
\mathcal{D}_k}\mathbf{P}_{\vartheta }
\{\mathit{dist}(\Gamma _{N},\vartheta ) >\epsilon \} \\[3mm]
&=\infl_{\displaystyle\vartheta _{N}\in
\mathcal{M}_{N}(\mathcal{D}^{\star})}
\supl_{\displaystyle\vartheta\in
\mathcal{D}_k}\mathbf{P}_{\vartheta } \{\mathit{dist}(\vartheta
_{N},\vartheta ) >\epsilon \} +\supl_{\displaystyle\vartheta\in
\mathcal{D}_k}\mathbf{P}_{\vartheta }
\{\mathit{dist}(\Gamma _{N},\vartheta ) >\epsilon \} \\[3mm]
&\le\infl_{\displaystyle\vartheta _{N}\in
\mathcal{M}_{N}(\mathcal{D}^{\star})}
\supl_{\displaystyle\vartheta\in
\mathcal{D}_k}\mathbf{P}_{\vartheta } \{\mathit{dist}(\vartheta
_{N},\vartheta ) >\epsilon \} + \varkappa
\end{array}
$$
So,
$$
\begin{array}{ll}
&\varkappa+\infl_{\displaystyle\vartheta _{N}\in
\mathcal{M}_{N}(\mathcal{D}^{\star})}
\supl_{\displaystyle\vartheta\in
\mathcal{D}_k}\mathbf{P}_{\vartheta }
\{\mathit{dist}(\vartheta _{N},\vartheta ) >\epsilon \}\ge\\[3mm]
&\ge\infl_{\displaystyle\vartheta _{N}\in
\mathcal{M}_{N}(\mathcal{D}_{k})} \supl_{\displaystyle\vartheta\in
\mathcal{D}_k}\mathbf{P}_{\vartheta }
\{\mathit{dist}(\vartheta _{N},\vartheta ) >\epsilon \}\ge\\[3mm]
&\infl_{\displaystyle\vartheta _{N}\in
\mathcal{M}_{N}(\mathcal{D}^{\star})}
\supl_{\displaystyle\vartheta\in
\mathcal{D}_k}\mathbf{P}_{\vartheta } \{\mathit{dist}(\vartheta
_{N},\vartheta ) >\epsilon \},
\end{array}
$$
and this is the fact we wanted to show.

By the definition of \emph{dist}, we have on the set
$\mathcal{M}_{N}(\mathcal{D}_{k})$:
$$
\mathit{dist}(\vartheta _{N},\vartheta ) = \| \vartheta _{N} -
\vartheta \| ^{(k)}.
$$

Further, for any $i=1,\ldots,k$ the following inclusion holds
$$
\{\| \vartheta _{N}-\vartheta \| ^{(k)} >\epsilon , \vartheta
_{N}\in \mathcal{D}_{k} \} \supseteq  \{|
\theta_{i}(N)-\theta_{i}| >\epsilon , \vartheta _{N}\in
\mathcal{D}_{k} \},
$$
where $\theta _{i}(N)$ is the \emph{i}-th component of the vector
$\vartheta _{N}$.

Therefore,
$$
\mathbf{P}_{\vartheta}\{\|\vartheta_{N}-\vartheta\|^{(k)}
>\epsilon, \vartheta _{N} \in \mathcal{D}_{k} \}
\ge \max_{1\le i\le k}\mathbf{P}_{\vartheta } \{|\theta
_{i}(N)-\theta _{i}|>\epsilon , \vartheta _{N}\in \mathcal{D}_{k}
\}.
$$

But estimation of the value
$$
\liminf_{N\to\infty}N^{-1} \ln \infl_{\displaystyle\vartheta
_{N}\in
\mathcal{M}_{N}(\mathcal{D}_{k})}\supl_{\displaystyle\vartheta\in
\mathcal{D}_{k}} \mathbf{P}_{\vartheta }\{| \theta _{i}(N)
-\vartheta ^{i}|  >\epsilon , \vartheta _{N}\in \mathcal{D}_{k} \}
 \stackrel{\triangle}{=}  A_{i}
$$
is exactly the problem already considered in the proof of Theorem
1 for the case of unique change-point. Therefore,
$$
A_{i}\ge -\min\left(
\intl_{\theta_i}^{\theta_i+\epsilon}\,J^{i-1}(t)dt,\,\intl_{\theta_i-\epsilon}^{\theta_i}\,J_{i}(t)dt\,
\right).
$$

So, finally we obtain
$$
\begin{array}{ll}
&\liminf_{N\to\infty}N^{-1}
\ln\infl_{\displaystyle\vartheta_{N}\in
\mathcal{M}_{N}(\mathcal{D}^{\star})}\supl_{\displaystyle\vartheta\in
\mathcal{D}_{k}}\mathbf{P}_{\vartheta }\{\| \vartheta
_{N}-\vartheta \| ^{(k)} >\epsilon , \vartheta _{N}\in
\mathcal{D}_{k}\} \ge\\[2mm]
&\ge -\minl_{1\le i\le k} \min\left(
\intl_{\theta_i}^{\theta_i+\epsilon}\,J^{i-1}(t)dt,\,
\intl_{\theta_i-\epsilon}^{\theta_i}\,J_{i}(t)dt\, \right).
\end{array}
$$
This completes the proof.

\section{Proof of Theorem 3} Due to the assumptions, the matrix
$I=\intl_0^1\,F(t)F^{*}(t)dt$ is positive definite. Therefore,
there exists the matrix $\big[N
\left(\mathcal{P}_1^N\right)^{-1}\big]$ for all $N>N_0(F)$. The
constant $N_0(F)$ can be exactly estimated for any given family of
functions $F(t)$.

Let us consider the matrix random process with continuous time
$\mathcal{Z}_{\scriptscriptstyle
N}(t)\df\mathcal{Z}_{\scriptscriptstyle N}([Nt]),\,t\in [0,1]$.

It is easy to see that the mathematical expectation of the process
$\mathcal{Z}_{\scriptscriptstyle N}(t)$ can be written as follows:
$$\begin{array}{ll}
&\mathbf{E}_{\theta}\mathcal{Z}_{\scriptscriptstyle N}(t)=N^{-1}\,
\left(\suml_{i=1}^{[Nt]}\,F(i/N)F^{*}(i/N)\,\mathbf{\Pi}^{*}(\theta,i)\right.\\
&\left.-\mathcal{P}_1^{[Nt]}(\mathcal{P}_1^N)^{-1}\,\suml_{i=1}^N\,
F(i/N)F^{*}(i/N)\,\mathbf{\Pi}^{*}(\theta,i)\right).
\end{array}
$$

After simple transformations we obtain that
$m(t)\df\liml_{N\to\infty}\mathbf{E}_{\theta}{\cal
Z}_{\scriptscriptstyle N}(t)$ has the form:
$$
m(t)=\left \{
\begin{array}{ll}
& A(t)I^{-1}(I-A(\theta))(\mathbf{a}-\mathbf{b})^{*}, \quad t\le\theta \\
& (I-A(t))I^{-1}A(\theta)(\mathbf{a}-\mathbf{b})^{*}, \quad t>
\theta,
\end{array}
\right. \eqno(C.0)
$$

Consider the square of the Gilbert norm of the matrix $m(t)$,
i.e., the function $f(t)=\textrm{tr} (m^{*}(t)m(t))$, and show
that the function $f(t)$ has a unique global maximum on the
segment $[0,1]$ at the point $t=\theta$.

First, for each $t\le\theta$:
$$
f(\theta)-f(t)=\textrm{tr} (B^{*}(A^2(\theta)-A^2(t))B),
$$
where matrix $B$ was defined in Theorem 3. Consider the matrix
$$
A^2(\theta)-A^2(t)=A(\theta)(A(\theta)-A(t))+(A(\theta)-A(t))A(t).
$$

Denote $L=A(\theta)(A(\theta)-A(t))$  and prove that the matrix
$L$ is positive definite as $t<\theta$. In fact, since the matrix
$A(\theta)$ is symmetric and positive definite, we can write
$$
x^{*}Lx=x^{*}A^{1/2}(\theta)A^{1/2}(\theta)(A(\theta)-A(t))x=
y^{*}\,A^{1/2}(\theta)(A(\theta)-A(t))A^{-1/2}(\theta)\,y,
$$
where $y=A^{1/2}(\theta)x$.

The matrices $A(\theta)-A(t)$ and
$A^{1/2}(\theta)(A(\theta)-A(t))A^{-1/2}(\theta)$ have identical
characteristic polynomial  and eigenvalues. Besides,
$A(\theta)-A(t)$ is positive definite as $t<\theta$. Therefore,
the matrix $A^{1/2}(\theta)(A(\theta)-A(t))A^{-1/2}(\theta)$ is
also positive definite as $t<\theta$ and therefore, the matrix $L$
is positive definite.

In analogy, the matrix $(A(\theta)-A(t))A(t)$ is positive definite
as $t<\theta$. Therefore, the matrix $A^2(\theta)-A^2(t)$ is
positive definite as $t<\theta$.

Now consider the matrix $D=B(A^2(\theta)-A^2(t))B^{*}$. The matrix
$D$ is positive definite if $\textrm{rank}(B)=M$, but this is our
assumption.

So, we obtain $\textrm{tr}(B(A^2(\theta)-A^2(t))B^*)>0$ for $t<
\theta$ and therefore, the function $f(t)$ has a unique global
maximum on the segment $[0,\theta]$ at the point $t=\theta$.

The same considerations for $t< \theta$ yield that $f(t)$
monotonically decreases on the segment $[\theta,1]$. As a result,
we obtain that $f(t)$ has a unique global maximum on the segment
$[0,1]$ at the point $t=\theta$.

Further, we are going to show the following: there exists a
positive constant $c$ such that $f(\theta)-f(t)\ge c \cdot |
\theta-t|$. This estimate can be obtained as follows. Taking into
account the continuity of the functions $f_j(t)$, we obtain
$$
A(\theta)-A(t)=\intl_t^{\theta}\,F(\tau)F^{*}(\tau)\,d\tau=
(\theta-t)U(t,\theta)>0, \eqno(C.1)
$$
where the matrix  $U(t,\theta)$ is positive definite for $0\le
t<\theta$ and negative definite for $t> \theta$. Due to the
continuity, we can write
$$
U(t,\theta)=U(\theta,\theta)+\kappa(t,\theta), \eqno(C.2)
$$
where $\kappa(t,\theta)\to 0$ as $t\to\theta$.

Then
$$\begin{array}{ll}
& f(\theta)-f(t)=\textrm{tr}
\left(B^*(A^2(\theta)-A^2(t))B\right)=\\[2mm]
& = \textrm{tr} \left(
BB^{*}A(\theta)(A(\theta)-A(t))\right)+\textrm{tr}
\left(BB^{*}(A(\theta)-A(t))A(t)\right)= \\[2mm]
& = (\theta-t)\,\textrm{tr}\left((\mathbf {a}-\mathbf
{b})^{*}(\mathbf {a}-\mathbf {b})V(t,\theta)\right),
\end{array}
\eqno(C.3)
$$
where
$V(t,\theta)=\left(E-A(\theta)I^{-1}\right)\left(A(\theta)U(t,\theta)+U(t,\theta)A(t)\right)
\left(E-I^{-1}A(\theta)\right)$.

Taking into account (C.1) and (C.2), we have
$$
\begin{array}{ll}
&V(t,\theta)=\left(E-A(\theta)I^{-1}\right)
\left(A(\theta)U(t,\theta)+U(t,\theta)A(t)\right)\left(E-I^{-1}A(\theta)\right)=\\[2mm]
&=\left(E-A(\theta)I^{-1}\right)
\left(A(\theta)U(\theta,\theta)+U(\theta,\theta)A(\theta)\right)
\left(E-I^{-1}A(\theta)\right)+\\[2mm]
&+\left(E-A(\theta)I^{-1}\right)
\left(A(\theta)\kappa(t,\theta)+\kappa(t,\theta)A(\theta)\right)
\left(E-I^{-1}A(\theta)\right)+\\[2mm]
&+(t-\theta)\left(E-A(\theta)I^{-1}\right)U(t,\theta)U(t,\theta)
\left(E-I^{-1}A(\theta)\right).
\end{array}
\eqno(C.4)
$$
Denote
$$\begin{array}{ll}
&
G(\theta)=\left(E-A(\theta)I^{-1}\right)\big(A(\theta)U(\theta,\theta)
+U(\theta,\theta)A(\theta)\big)
\left(E-I^{-1}A(\theta)\right)\\[2mm]
& R(t,\theta)=\left(E-A(\theta)I^{-1}\right)
\big(A(\theta)\kappa(t,\theta)+\kappa(t,\theta)A(\theta)\big)
\left(E-I^{-1}A(\theta)\right)\\[2mm]
& H(t,\theta)=\left(E-A(\theta)I^{-1}\right)U(t,\theta)U(t,\theta)
\left(E-I^{-1}A(\theta)\right)
\end{array}
\eqno(C.5)
$$
and put
$$
\tilde{G}(\theta)=\left \{
\begin{array}{ll}
& G(\theta),\quad \theta>t \\
& -G(\theta), \quad \theta\le t.
\end{array}
\right. \eqno(C.6)
$$

Then from (C.3), (C.4), (C.5) and (C.6) we get
$$
\begin{array}{ll}
& f(\theta)-f(t) =|\theta-t|\,\textrm{tr}\left((\mathbf
{a}-\mathbf {b})^*(\mathbf {a}-\mathbf
{b})\tilde{G}(\theta)\right)+\\[2mm]
&+(\theta-t)\,\textrm{tr}\left((\mathbf {a}-\mathbf {b})^*(\mathbf
{a}-\mathbf {b})R(t,\theta)\right)-\\[2mm]
& -(\theta-t)^2\,\textrm{tr}\left((\mathbf {a}-\mathbf
{b})^*(\mathbf {a}-\mathbf {b})H(t,\theta)\right)
\end{array}
\eqno(C.7)
$$

Since $R(t,\theta)\to 0$ as $t\to\theta$ and $H(t,\theta)$ is
positive definite, we conclude that
$$
f(\theta)-f(t)\ge |\theta-t|\,\textrm{tr}\left((\mathbf
{a}-\mathbf {b})^*(\mathbf {a}-\mathbf
{b})\tilde{G}(\theta)\right)+o(|t-\theta|),
$$
i.e., there exists a positive definite matrix $W(\theta)$ such
that
$$
\|m(\theta)\|^2-\|m(t)\|^2=f(\theta)-f(t)\ge
|\theta-t|\,\textrm{tr}\big((\mathbf {a}-\mathbf {b})^*(\mathbf
{a}-\mathbf {b})W(\theta)\big)
$$
for some neighborhood of $\theta$.  Therefore, we have got the
estimate of sharpness of the maximum for the function $f(t)$:
$$
f(\theta)-f(t)\ge |\theta-t|\lambda_{\scriptscriptstyle
F}\mathrm{tr}\big[(\mathbf {a}-\mathbf {b})^*(\mathbf {a}-\mathbf
{b})\big], \eqno(C.8)
$$
where
$$
\lambda_{\scriptscriptstyle
F}\df\minl_{\beta\le\theta\le\alpha}\dpfrac{\textrm{tr}\big[(\mathbf
{a}-\mathbf {b})^*(\mathbf {a}-\mathbf
{b})W(\theta)\big]}{\textrm{tr}\big[(\mathbf {a}-\mathbf
{b})^*(\mathbf {a}-\mathbf {b})\big]}.
$$

Let us describe how to calculate $\lambda_{\scriptscriptstyle F}$.
For given family of functions $F(t)$ we can calculate the function
$f(t)=\mathrm{tr}\big[m^*(t)m(t)\big]$. Then it is possible to
calculate
$$
\lambda_{\scriptscriptstyle F}=\minl_{\beta\le
t\le\alpha,\,\beta\le\theta\le\alpha}\dpfrac{f(\theta)-f(t)}
{|\theta-t|\textrm{tr}\big[(\mathbf {a}-\mathbf {b})^*(\mathbf
{a}-\mathbf {b})\big]}.
$$
Due to the condition $0<\beta\le\theta\le\alpha<1$, we get
$\lambda_{\scriptscriptstyle F}>0$ (see (C.5)). Note that from
(C.8) and definition of $f(t)$ we have for any $t\in
[\beta,\alpha]$:
$$
\|m(\theta)\|^2-\|m(t)\|^2 \ge\dpfrac{\lambda_{\scriptscriptstyle
F}}{2\|m(\theta\|}|\theta-t|\,\|\mathbf {a}-\mathbf {b}\|^2
\eqno(C.9)
$$

\bigskip
The process ${\cal Z}_{\scriptscriptstyle N}(t)$ can be decomposed
into deterministic and stochastic terms:
$$
{\cal Z}_{\scriptscriptstyle N}(t)=m(t)+\gamma_{\scriptscriptstyle
N}(t)+\eta_{\scriptscriptstyle N}(t), \eqno(C.10)
$$
where the norm of the deterministic function
$\gamma_{\scriptscriptstyle N}(t)$
 converges to zero with the rate $L_{F}/N)$ (this term estimates
the difference between corresponding integral sum and the
integral; the constant $L_{F}$ depends of the function family
$F(t)$ and \emph{can be estimated explicitly for any given
family}), and the stochastic term is equal to
$$
\eta_{\scriptscriptstyle
N}(t)=N^{-1}\left(\suml_{i=1}^{[Nt]}\,F(i/N)\nu^{*}_i- \mathcal
{P}_1^{[Nt]}(\mathcal
{P}_1^N)^{-1}\,\suml_{i=1}^N\,F(i/N)\nu^*_i\right).
$$

The norm of the process $\eta_{\scriptscriptstyle N}(t)$ can be
estimated as follows:
$$
\begin{array}{ll}
&\supl_{\beta\le t\le\alpha}\|\eta_{\scriptscriptstyle N}(t)\|\le
R
\left[\sqrt{K}+\|I\|\cdot\|I^{-1}\|+\frac{L_{\scriptscriptstyle F}}{N}\big(\|I\|+\|I^{-1}\|+L_{\scriptscriptstyle F}/N\big)\right]\times\\[2mm]
&\times\left(\maxl_{1\le i\le K}\maxl_{1\le l\le M}\maxl_{[\beta
N]\le
n\le N}\,N^{-1}|\suml_{j=1}^n\,f_i(j/N)\nu_{lj}|\right)\df\\[2mm]
&=\mathcal{R}\left(\maxl_{1\le i\le K}\maxl_{1\le l\le
M}\maxl_{[\beta N]\le n\le
N}\,N^{-1}|\suml_{j=1}^n\,f_i(j/N)\nu_{lj}|\right),
\end{array}
\eqno(C.11)
$$
where $\mathcal{R}=\mathcal{R}(F,N)$. Here we used the following
relations
$$
\begin{array}{ll}
&\maxl_{t\in [0,1]}\|N^{-1}\mathcal{P}_1^{[Nt]}-A(t)\|\le\frac{L_{\scriptscriptstyle F}}{N},\,\,\maxl_{t\in [0,1]}\|A(t)\|\le\|I\|\\[2mm]
 &\|N(\mathcal{P}_1^N)^{-1}-I^{-1}\||\le\frac{L_{\scriptscriptstyle F}}{N}
\end{array}
$$
and took into account that for any matrix $M$ we have the relation
$\|M\|=\sqrt{\mathrm{tr}(M^*M)}\le R\maxl_{i,j}|m_{ij}|$, where
constant $R$ depends only of the dimensionality.

Denote $\tilde{S}_n=\suml_{j=1}^n\,f_i(j/N)\nu_{lj},\,\,\tilde{\xi}(j)=f_i(j/N)\nu_{lj}$ and\\
 put $\sigma^2=\supl_i\supl_{1\le n\le N}\supl_{1\le l\le
M}\mathbf{E}_{\theta}(f_i(n/N)\nu_{ln})^2$. Choose the number
$\epsilon(x)$ from the following condition
$$
\ln(1+\epsilon(x))=\left \{
\begin{array}{ll}
& x^2/4g,\quad x\le g T, \\
& xT/4, \quad x>g T,
\end{array}
\right.
$$
where the constant $T$ is taken from the uniform Cramer condition
and $g>\sigma^2$.

For the chosen $\epsilon(x)=\epsilon$, we choose the number
$m_0(x)\ge 1$ from the uniform $\psi$-mixing condition such that
$\psi(m)\le \epsilon$ for $m\ge m_0(x)$.

Decompose the sum  $\tilde S_n$ into groups of weakly dependent
terms:
$$\tilde
S_n=\tilde S_n^1+\tilde S_n^2+\dots+\tilde S_n^{m_0(x)},
$$
where
$$
\tilde S_n^i=\tilde\xi(i)+\tilde \xi(i+m_0(x))+\dots+\tilde\xi
\left(i+m_0(x)[\frac {n-i}{m_0(x)}]\right),
$$
and $i=1,2,\dots,m_0(x)$.

The number of summands  $k(i)$ in each group is no less than
$[n/m_0(x)]$ and no more than $[n/m_0(x)]+1$. The $\psi$-mixing
coefficient between summands within each group is no larger than
$\epsilon$. Therefore,
$$
\begin{array}{ll}
& \mathbf{P}_{\theta}\{|\tilde{S}_n|/n \ge x\}\le
\suml_{i=1}^{m_0(x)}\,\mathbf{P}_{\theta}\{|\tilde
{S}_n^i/n| \ge x/m_0(x)\}\le \\
& \le m_0(x) \,\maxl_{1\le i\le
m_0(x)}\mathbf{P}_{\theta}\{|\tilde {S}_n^i| \ge (k(i)-1)x\}.
\end{array}
\eqno(C.12)
$$

From Chebyshev's inequality we have:
$$
\mathbf{P}_{\theta}\left\{\tilde{S}^i_k=\suml_{j=0}^k\tilde{\xi}(i+m_0j)\ge
x\right\} \le e^{-tx}\mathbf{E}_{\theta}e^{t\tilde{S}^i_k},\quad
\forall t>0. \eqno(C.13)
$$

Further, from  $\psi$-mixing condition it follows that (see
Ibragimov, Linnik (1971)):
$$
\mathbf{E}_{\theta}e^{t\tilde{S}^i_k}\le
(1+\epsilon)^k\,\mathbf{E}_{\theta}\exp(t\tilde\xi(i))\mathbf{E}_{\theta}\exp(t\tilde\xi(i+m_0))\dots
\mathbf{E}_{\theta}\exp(t\tilde\xi(i+m_0 k)). \eqno(C.14)
$$

Consider the term $\mathbf{E}_{\theta}\exp(t\tilde\xi(i))$. From
the uniform Cramer's condition it follows that for each $0<t< T$:
$$
\mathbf{E}_{\theta}e^{t\tilde \xi(i)}\le \exp (t^2g/2).
$$
Then from (C.13) and (C.14) we obtain
$$
\mathbf{P}_{\theta}\{\tilde{S}_k^i \ge x\}\le (1+\epsilon)^k \exp
\left(kgt^2/2-tx\right).
$$

Taking the minimum of $kgt^2/2-tx$ w.r.t. $t$, write
$$
\mathbf{P}_{\theta}\{\tilde{S}_k^i\ge x\}\le \left \{
\begin{array}{ll}
(1+\epsilon)^k\,\exp(-x^2/2kg), \quad x\le kgT, \\
(1+\epsilon)^k\,\exp(- xT/2), \quad x>kgT.
\end{array}
\right.
$$

From the definition of  $\epsilon$ we obtain
$$
\mathbf{P}_{\theta}\{|\tilde{S}_k^i/k|\ge x\}\le \left \{
\begin{array}{ll}
\exp(-kx^2/4g), \quad x\le gT, \\
\exp(-kxT/4), \quad x>gT.
\end{array}
\right. \eqno(C.15)
$$

Now, using (C.12) and (C.15), we obtain
$$
\mathbf{P}_{\theta}\{|\tilde S_n/n|\ge x\}\le \left \{
\begin{array}{ll}
& m_0(x)\exp\left(-x^2 n/4gm_0(x)\right), \quad x\le gT, \\
& m_0(x) \exp \left(-Txn/4m_0(x)\right), \quad x>gT.
\end{array}
\right. \eqno(C.16)
$$

From (C.11) and (C.16) we get
$$
\mathbf{P}_{\theta}\{\supl_{\beta\le
t\le\alpha}\|\eta_{\scriptscriptstyle N}(t)\|> \epsilon \} \le
m_0(\epsilon/\mathcal{R}) \left \{
\begin{array}{ll}
& \exp \left(-(\epsilon/\mathcal{R})^2 N\beta/4g m_0(\epsilon/\mathcal{R})\right),\\
& \epsilon\le \mathcal{R}gT \\
& \exp\left(-T(\epsilon/\mathcal{R})
N\beta/4m_0(\epsilon/\mathcal{R})\right), \\
& \epsilon >\mathcal{R}gT,
\end{array}
\right. \eqno(C.17)
$$
In particular, for the case of independent observations,
$m_0(\epsilon)=1$.

From the definition of the estimate
$\hat{\theta}_{\scriptscriptstyle N}$ and (C.9) we can write
$$
\begin{array}{ll}
& \mathbf{P}_{\theta}\left\{|\hat{\theta}_{\scriptscriptstyle
N}-\theta|> \epsilon,\, \hat{\theta}_{\scriptscriptstyle
N}\in\mathrm{Arg}\maxl_{\beta\le t\le \alpha
}\|\mathcal{Z}_{\scriptscriptstyle N}(t)\|
\right\}=\\[2mm]
&=\mathbf{P}_{\theta}\{\|\mathcal{Z}_{\scriptscriptstyle
N}(\hat\theta_{\scriptscriptstyle N})\|\ge
\|\mathcal{Z}_{\scriptscriptstyle N}(t)\|,\,t\in[\beta,\alpha],
\,|\hat\theta_{\scriptscriptstyle N}-\theta| > \epsilon \} \\[2mm]
& \le \mathbf{P}_{\theta} \{\|\eta_{\scriptscriptstyle
N}(\hat{\theta}_{\scriptscriptstyle N})\|-
\|\eta_{\scriptscriptstyle N}(\theta)\|\ge
\|m(\theta)\|^2-\|m(\hat{\theta}_{\scriptscriptstyle
N})\|^2+L_{\scriptscriptstyle F}/N,\;
|\hat\theta_{\scriptscriptstyle N}-\theta| > \epsilon \} \\[2mm]
& \le \mathbf{P}_{\theta}\left\{\supl_{\beta\le t\le\alpha}
\|\eta_{\scriptscriptstyle N}(t)\| \ge
\Big[\dpfrac{\epsilon\lambda_{\scriptscriptstyle
F}}{4\|m(\theta\|}\textrm
{tr}((\mathbf{a}-\mathbf{b})^{*}(\mathbf{a}-\mathbf{b}))-\dpfrac{L_{\scriptscriptstyle
F}}{N}\Big] \right\}\le\\[2mm]
&\le\mathbf{P}_{\theta}\left\{\supl_{\beta\le t\le\alpha}
\|\eta_{\scriptscriptstyle N}(t)\| \ge
\Big[\dpfrac{\epsilon\lambda_{\scriptscriptstyle
F}}{4\mathcal{M}}\textrm
{tr}((\mathbf{a}-\mathbf{b})^{*}(\mathbf{a}-\mathbf{b}))-\dpfrac{L_{\scriptscriptstyle
F}}{N}\Big] \right\},
\end{array}
\eqno(C.18)
$$
where $\mathcal{M}=\maxl_{\beta\le\theta\le\alpha}\|m(\theta\|$.

 Denote $C(\epsilon,N)=\Big[\dpfrac{\epsilon\lambda_{\scriptscriptstyle F}}{4\mathcal{M}}
\|\mathbf{a}-\mathbf{b}\|^2-\dpfrac{L_{\scriptscriptstyle
F}}{N}\Big]$. Then, finally we obtain from (C.18):
$$
\supl_{\beta\le\theta\le\alpha}\mathbf{P}_{\theta}\{|\hat\theta_{\scriptscriptstyle
N}-\theta|> \epsilon \} \le
m_0\left(C(\epsilon,N)/\mathcal{R}\right)\left \{
\begin{array}{ll}
& \exp \left(-\dpfrac {N\beta
\Big(C(\epsilon,N)/\mathcal{R}\Big)^2}
{4gm_0\left(C(\epsilon,N)/\mathcal{R}\right)}\right),\\
& \textrm{  if  } C(\epsilon,N) \le \mathcal{R}gT \\
& \exp\left(-\dpfrac
{TN\beta\Big(C(\epsilon,N)/\mathcal{R}\Big)}{4m_0\left(C(\epsilon,N)/\mathcal{R}\right)}\right),\\
&\textrm{ if } C(\epsilon,N)>\mathcal{R}gT.
\end{array}
\right.
$$
\begin{rmk}
In case of only one regression relationship and independent noises
$\nu_i$, we obtain from here
$$
\supl_{\beta\le\theta\le\alpha}\mathbf{P}_{\theta}\{|\hat{\theta}_{\scriptscriptstyle
N}-\theta|> \epsilon \} \le \left \{
\begin{array}{ll}
& \exp \left(-\dpfrac {N\beta\epsilon^2}{4g\mathcal{R}^2}
\Big[\dpfrac{\lambda_{\scriptscriptstyle F}}{4\mathcal{M}}
\suml_{j=1}^k\,(a_j-b_j)^2)-\dpfrac{L_{\scriptscriptstyle F}}{N}
\Big]^2\right)\\
& \textrm{if}\,\, C(\epsilon,N) \le \mathcal{R}gT \\
&
\exp\left(-\dpfrac{TN\beta\epsilon}{4\mathcal{R}}\Big[\dpfrac{\lambda_{\scriptscriptstyle
F}}{4\mathcal{M}}
\suml_{j=1}^k\,(a_j-b_j)^2)^2-\dpfrac{L_{\scriptscriptstyle F}}{N}
\Big]\right) \\
&\textrm{ if }\,\, C(\epsilon,N)>\mathcal{R}gT.
\end{array}
\right.
$$
\end{rmk}

Theorem 3 is proved.

\medskip

Corollary 2 can be obtained (as it follows from the proof) from
the estimates of $\mathbf{P}_{\theta}\{\supl_{\beta\le t\le\alpha}
\|\eta_{\scriptscriptstyle
N}(t)\|>\epsilon\},\,\theta=0\,\,\mbox{or}\,\,\theta\not=0$.

\section{ Proof of Theorem 4}
The proof is based on the same ideas as in Section C, and so we
give the sketch of the proof.

Let us consider the matrix random process with continuous time
$\mathbb{Z}_{\scriptscriptstyle N}(t)
\df\mathbb{Z}_{\scriptscriptstyle N}([Nt]),\,t\in [0,1]$.

It is easy to see that the mathematical expectation of the process
$\mathbb{Z}_{\scriptscriptstyle N}(t)$ can be written as follows:
$$
\mathbf{E}_{\theta}\mathbb{Z}_{\scriptscriptstyle
N}(t)=N^{-1}\left(\suml_{n=1}^{[Nt]}V(n/N)\mathbf{\Pi}^*(\theta,n)-\mathcal{T}_{1}^{[Nt]}
(\mathcal{T}_1^N)^{-1}\suml_{n=1}^{N}V(n/N)\mathbf{\Pi}^*(\theta,n)\right)
$$

Denote
$M(t)\df\liml_{N\to\infty}\mathbf{E}_{\theta}\mathbb{Z}_{\scriptscriptstyle
N}(t)$. After simple transformation we have
$$
M(t)=\left\{\begin{array}{ll}
&\mathbb{R}(t)\mathbb{R}^{-1}\left(\mathbb{R}-\mathbb{R}(\theta)\right)(\mathbf{a}-\mathbf{b})^*,\,\,t\le\theta\\[2mm]
&\left(\mathbb{R}-\mathbb{R}(t)\right)\mathbb{R}^{-1}\mathbb{R}(\theta)(\mathbf{a}-\mathbf{b})^*,\,\,t>\theta
\end{array}
\right. \eqno(D.1)
$$
It can be shown from (D.1) (by the analogous arguments as in
Section C) that the function $\Phi(t)\df\|M(t)\|^2=
\mathrm{tr}\left(M(t)M^*(t)\right)$ has unique global maximum on
the segment $[0,1]$ at the point $t=\theta$ and there exists
$\lambda_{\scriptscriptstyle V}>0$ such that the following
inequality holds
$$
\Phi(\theta)-\Phi(t)\ge \lambda_{\scriptscriptstyle
V}|\theta-t|\mathrm{tr}\left[(\mathbf{a}-\mathbf{b})(\mathbf{a}-\mathbf{b})^*\right]
\eqno(D.2)
$$
for any $\beta\le t\le\alpha$. The constant
$\lambda_{\scriptscriptstyle V}$ depends only of $V(t)$ and can be
estimated analogously the constant $\lambda_{\scriptscriptstyle
F}$ from Section C.

Consider matrix sequence $N^{-1}\mathcal{T}^{\scriptscriptstyle
N}_1$. Due to the assumptions, this sequence
$\mathbf{P}_{\theta}$-a.s. tends to the positive definite matrix
$\mathbb{R}=\intl_0^1V(s)ds$, and the rate of the convergence is
exponential. Therefore, there exists number
$N_1=N_1\left(\{\mathbf{X}(n)\}\right)$ such that as $N>N_1$ we
get
$$
\mathbf{P}_{\theta}\{\|N^{-1}\mathcal{T}^{\scriptscriptstyle
N}_1-\mathbb{R}\|>\epsilon\}\le
L(\epsilon)\exp\left(-K(\epsilon)N\right), \eqno(D.3)
$$
where functions $L(\epsilon),\,K(\epsilon)$ can be exactly
estimated (taking into account $\psi$-mixing condition and
Cramer's condition) by the scheme of Section C. The number $N_1$
can be estimated by the random sequence $\{\mathbf{X}(n)\}$.

Process $\mathbb{Z}_{\scriptscriptstyle N}(t)$ can be written as
follows
$$
\mathbb{Z}_{\scriptscriptstyle
N}(t)=M(t)+\Gamma_{\scriptscriptstyle
N}(t)+\zeta_{\scriptscriptstyle N}(t),
$$
where $\Gamma_{\scriptscriptstyle
N}(t)=\mathbf{E}_{\theta}\mathbb{Z}_{\scriptscriptstyle
N}(t)-M(t)$ and $\zeta_{\scriptscriptstyle
N}=\mathbb{Z}_{\scriptscriptstyle
N}(t)-\mathbf{E}_{\theta}\mathbb{Z}_{\scriptscriptstyle N}(t)$.

Note that $\maxl_{0\le t\le 1}\|\Gamma_{\scriptscriptstyle
N}(t)\|\le\frac{L_{\scriptscriptstyle V}}{N}$ (because this is the
difference between the sum and the integral), and constant
$L_{\scriptscriptstyle V}$ can be estimated exactly for any given
function $V(t)$.

Fix $\epsilon,\,\,0<\epsilon<\min\left((\alpha-\beta),
\|\mathbb{R}\|/2\right)$ and consider the events
$$
\begin{array}{ll}
&D_{\scriptscriptstyle N}=\{
\|N^{-1}\mathcal{T}^{\scriptscriptstyle
N}_1-\mathbb{R}\|\le\|\mathbb{R}\|/2,
\\[3mm]
&\maxl_{0\le t\le 1}
\|N^{-1}\mathcal{T}_1^{[Nt]}-\mathbb{R}(t)\|<\epsilon,\,\,\|N(\mathcal{T}_1^N)^{-1}
-\mathbb{R}^{-1}\|<
\epsilon \},\\[3mm]
&\bar{D}_{\scriptscriptstyle N}=\Omega\backslash
D_{\scriptscriptstyle N}.
\end{array}
$$
Note that matrix $N^{-1}\mathcal{T}_1^N$ is non-degenerate on the
set $D_{\scriptscriptstyle N}$. Then, due to (D.3),
$$
\delta_{\scriptscriptstyle
N}(\epsilon)\df\mathbf{P}_{\theta}(\bar{D}_{\scriptscriptstyle
N})\le 3L(\epsilon)\exp\left(-K(\epsilon)N\right). \eqno(D.4)
$$

Further, analogously (C.11), we can write on the set
$D_{\scriptscriptstyle N}$
$$
\begin{array}{ll}
&\supl_{\beta\le t\le\alpha}\|\zeta_{\scriptscriptstyle N}(t)\|\le
R \left[\sqrt{K}+\|\mathbb{R}\|\cdot\|\mathbb{R}^{-1}\|+\epsilon
\big(\|\mathbb{R}\|+\|\mathbb{R}^{-1}\|+\epsilon\big)\right]\times\\[2mm]
&\times\left(\maxl_{1\le i\le K}\maxl_{1\le l\le M}\maxl_{[\beta
N]\le
n\le N}\,N^{-1}|\suml_{j=1}^n\,x_{ij}\nu_{lj}|\right)\df\\[2mm]
&=\mathbf{R}\left(\maxl_{1\le i\le K}\maxl_{1\le l\le
M}\maxl_{[\beta N]\le n\le
N}\,N^{-1}|\suml_{j=1}^n\,x_{ij}\nu_{lj}|\right),
\end{array}
\eqno(D.5)
$$
where $\mathbf{R}=\mathbf{R}(V,\epsilon)$.

Now we can use (C.17) and get (by the analogous reasons) from
(D.5) on the set $D_{\scriptscriptstyle N}$
$$
\mathbf{P}_{\theta}\{\supl_{\beta\le
t\le\alpha}\|\zeta_{\scriptscriptstyle N}(t)\|>
\epsilon,\,\,\mathbb{I}(D_{\scriptscriptstyle N} \} \le
m_0(\epsilon/\mathbf{R}) \left \{
\begin{array}{ll}
& \exp \left(-(\epsilon/\mathbf{R})^2 N\beta/4g m_0(\epsilon/\mathbf{R})\right),\\
& \epsilon\le \mathbf{R}gT \\
& \exp\left(-T(\epsilon/\mathbf{R})
N\beta/4m_0(\epsilon/\mathbf{R})\right),\\
& \epsilon > \mathbf{R}gT,
\end{array}
\right. \eqno(D.6)
$$

Using  (D.4), (D.6), and the analogous considerations as in
(C.18), we get
$$
\begin{array}{ll}
&\supl_{\beta\le\theta\le\alpha}\mathbf{P}_{\theta}\{|\hat\theta_{\scriptscriptstyle
N}-\theta|> \epsilon
\} \le \delta_{\scriptscriptstyle N}(\epsilon)+\\[3mm]
&m_0\left(\mathbb{C}(\epsilon,N)/\mathbf{R}\right)\left \{
\begin{array}{ll}
& \exp \left(-\dpfrac {N\beta
\Big(\mathbb{C}(\epsilon,N)/\mathbf{R}\Big)^2}
{4gm_0\left(\mathbb{C}(\epsilon,N)/\mathbf{R}\right)}\right), \textrm{  if  } \mathbb{C}(\epsilon,N) \le \mathbf{R}gT \\
& \exp\left(-\dpfrac
{TN\beta\Big(\mathbb{C}(\epsilon,N)/\mathbf{R}\Big)}{4m_0\left(\mathbb{C}(\epsilon,N)/\mathbf{R}\right)}\right),
\textrm{ if } \mathbb{C}(\epsilon,N)>\mathbf{R}gT,
\end{array}
\right.
\end{array}
$$
where
$\mathbb{C}(\epsilon,N)=\Big[\dpfrac{\epsilon\lambda_{\scriptscriptstyle
V}}{4\mathbb{M}}
\|\mathbf{a}-\mathbf{b}\|^2-\dpfrac{L_{\scriptscriptstyle
V}}{N}\Big],\,\mathbb{M}=\maxl_{\beta\le t\le\alpha}\|M(t)\|$.

Theorem 4 is proved.
\section{Proof of Theorem 5}

The proposed method of multiple change-point detection and
estimation is based upon the idea of recurrent reduction to the
case of one change-point.

In order to prove theorem 5 we need to prove the following two
propositions:

i) in the case of a stationary sub-sample the norm of the decision
statistic does not exceed the threshold with the great
probability. This fact is exactly the result of Corollary 2;

ii) in the case of a non-stationary sub-sample with at least two
change-points, the norm of the decision statistic exceeds the
decision threshold with the great probability.

In order to illustrate ii), let us consider a sub-sample of size
$N$ with two change-points $0<\theta_1<\theta_2 <1$.

In this case the decision statistic can be decomposed into a
deterministic and a stochastic term (see (C.10)).

We have from (C.0) for $0\le t\le \theta_1$:
$$\begin{array}{ll}
& m(t)=A(t)a_1-A(t)A^{-1}(1)\,\big(A(\theta_1)a_1+A(\theta_1,\theta_2)a_2+A(\theta_2,1)a_{3}\big)\\[2mm]
&=A(t)\,(a_1-A^{-1}(1)u),
\end{array}
\eqno(E.1)
$$
where
$u=A(\theta_1)a_1+A(\theta_1,\theta_2)a_2+A(\theta_2,1)a_{3}$.

Again using (C.0), we get for $\theta_1\le t\le\theta_2$:
$$\begin{array}{ll}
&m(t)=A(\theta_1)a_1+A(\theta_1,t)a_2-A(t)A^{-1}(1)u=\\[2mm]
&=A(\theta_1)\big(a_1-A^{-1}(1)u\big)+A(\theta_1,t)\big(a_2-A^{-1}(1)u\big).
\end{array}
$$

If
$$
\|m(\theta_1)\| \ge \Lambda\df\dpfrac B{2(h+1)}>0,
$$
then $\maxl_{\beta\le t\le \alpha}\,\|m(t)\|\ge \Lambda >0$.

Otherwise, let $\|m(\theta_1)\| < \Lambda$. Then
$$\begin{array}{ll}
& \|m(\theta_2)\| \ge
\|A(\theta_1,\theta_2)(a_2-A^{-1}(1)(u)\|-\Lambda=\\[2mm]
&=\|A(\theta_1,\theta_2)(a_2-a_1+a_1-A^{-1}(1)u\|-\Lambda \\[2mm]
& \ge
\|A(\theta_1,\theta_2)(a_2-a_1)\|-\|A(\theta_1,\theta_2)(a_1-A^{-1}(1)(u)\|-\Lambda
\\[2mm]
&\ge
B-\|A(\theta_1,\theta_2)A^{-1}(\theta_1)\|\,\Lambda-\Lambda\ge
B-\Lambda(1+h)> \Lambda.
\end{array}
$$

Therefore, taking into account (E.1), we get: there exists
$\Lambda>0$ such that
$$
\maxl_{\beta\le t\le\alpha}\|m(t)\|\ge\Lambda \eqno(E.2)
$$
From (E.2) it follows that we get ii) with the great probability.

\medskip

After these preliminary considerations, let us consider the
probability of the event:
$$
(\hat k_N \ne k)\cup \{(\hat k_N = k)\cap (\maxl_{1\le i\le
k}\,|\hat\theta_{Ni}-\theta_i|> \delta) \eqno(E.3)
$$
for some fixed $\delta,\,\,\epsilon>\delta>0$. Let us consider the
following cases:

$a) \,\,\{\hat k_N < k\},\,\, \,  b)\,\,\{\hat k_N >k\},\,\,\,
c)\,\,\{(\hat k_N = k)\cap (\maxl_{1\le i\le
k}\,|\hat\theta_{Ni}-\theta_i|> \delta)\}$.

\medskip
\emph{Case a)}

In this case the proposed method does not detect at least one
change-point, i.e., a certain sub-sample of size $\tilde N\ge
[2\delta N]$ containing at least one true change-point, is
classified as stationary. Then
$$
\mathbf{P}_{\vartheta}\{\hat k_N <k\}\le
\mathbf{P}_{\vartheta}\{\maxl_{\beta \le t\le
\alpha}\,\|\mathcal{Z}_{\scriptscriptstyle \tilde{N}}(t)\|\le
C(\tilde N)\}
\eqno(E.4)
$$
where $C(\tilde{N})$ is the decision threshold for the sub-sample.

Choose  $C(\tilde{N})< \Lambda$. Then due to (E.4) and (C.10) we
have
$$
\begin{array}{ll}
&\mathbf{P}_{\vartheta}\{\maxl_{\beta\le
t\le\alpha}\|\mathcal{Z}_{\scriptscriptstyle \tilde{N}}\|\le
C(\tilde{N})\}\le
\mathbf{P}_{\vartheta}\{\maxl_{\beta\le t\le\alpha}\|\eta_{\scriptscriptstyle \tilde{N}}(t)\|\ge\maxl_{\beta\le t\le\alpha}\|m(t)\|-\frac{L_{\scriptscriptstyle F}}{N}-C(\tilde{N})\}\\[3mm]
&\le\mathbf{P}_{\vartheta}\{\maxl_{\beta\le
t\le\alpha}\|\eta_{\scriptscriptstyle
\tilde{N}}(t)\|\ge\Lambda-\frac{L_{\scriptscriptstyle
F}}{N}-C(\tilde{N})\}
\end{array}
$$
Now we can use (C.17), changing $\epsilon$ by
$\{\Lambda-\frac{L_{\scriptscriptstyle F}}{N}-C(\tilde{N})\}$, and
get the exponential estimate for the event
$\{\hat{k}_{\scriptscriptstyle N}<k\}$.

\medskip

\emph{Case b)}

In this case there exists a stationary sub-sample of the size
$\hat{N}\ge[\delta N]$ such that it is classified as
non-stationary. Then
$$
 \mathbf{P}_{0}\{\hat{k}_{\scriptscriptstyle N}>k\}\le\mathbf{P}_{0}
 \{\maxl_{\beta\le t\le\alpha}\|\mathcal{Z}_{\scriptscriptstyle \hat{N}}(t)\|>C(\hat{N})\}
 \eqno(E.5)
 $$
But the exponential estimate of the right-hand side (E.5) can be
taken from (9).

\medskip

\emph{Case c)}

In this case there exists a sub-sample of the size $N^*\ge
[2\delta N]$ such that the distance between a true change-point
parameter $\theta_i$ and its estimate $\hat\theta_{Ni}$ is larger
than $\delta$. This is exactly the case of Theorem 3, and we get
the exponential estimate of this event from (8).

\bigskip
Therefore, we get the exponential estimate for the event (E.3).
This completes the proof of Theorem 5.

\end{document}